\documentclass[11pt]{amsart}
\usepackage{amsmath,mathtools}
\usepackage{amsmath, amsthm, amssymb, amsfonts, enumerate}
\usepackage[colorlinks=true,linkcolor=blue,urlcolor=blue]{hyperref}
\usepackage{dsfont}
\usepackage{color}
\usepackage{geometry}
\usepackage{todonotes}
\usepackage{epstopdf}
\usepackage{bbm}
\usepackage{geometry}
\usepackage[utf8]{inputenc}
\usepackage[format=plain, indention=1cm]{caption}
\usepackage{amsfonts}
\usepackage{amsfonts}
\usepackage{textcomp}
\usepackage{pgfplots}
\usepackage{amssymb}
\usepackage{float}
\usepackage{tikz}
\usepackage{epsfig}
\usepackage{amsmath}
\usepackage[english]{babel}
\usepackage{a4}
\usepackage{enumerate}
\usepackage{amsmath}
\newcommand{\stkout}[1]{\ifmmode\text{\sout{\ensuremath{#1}}}\else\sout{#1}\fi}

\usepackage{setspace}
\usepackage{tabularx}

\newtheorem{theorem}{Theorem}[section]
\newtheorem{remark}[theorem]{Remark}
\newtheorem{assumption}[theorem]{Assumption}
\newtheorem{lemma}[theorem]{Lemma}
\newtheorem{proposition}[theorem]{Proposition}
\newtheorem{corollary}[theorem]{Corollary}

\def \E{\mathbb{E}}

\def \P{\mathsf{P}}

\def \R{\mathbb{R}}

\definecolor{red}{rgb}{1.0,0.0,0.0}

\definecolor{blu}{rgb}{0.0,0.0,1.0}

\definecolor{gre}{rgb}{0.03,0.50,0.03}

\renewcommand{\hat}{\widehat}
\newcommand{\eins}{\text{$\mathbbm{1}$}}
\renewcommand{\bar}{\overline}
\renewcommand{\phi}{\varphi}

\newcommand{\cvd}{$\quad\Box $
                  \medskip

                 }

\title[]{Uncertainty over Uncertainty in Environmental \\ Policy Adoption: 
Bayesian Learning of Unpredictable Socioeconomic Costs} 
\author[]{Matteo Basei$^*$}
\thanks{$^*$EDF R\&D Paris \& FIME (Finance for Energy Markets Research Centre), Paris, France. Email: matteo.basei@edf.fr}
\author[]{Giorgio Ferrari$^\dagger$}
\thanks{$^\dagger$Center for Mathematical Economics (IMW), Bielefeld University, Universit\"atsstrasse 25, 33615, Bielefeld, Germany. Email: giorgio.ferrari@uni-bielefeld.de}
\author[]{Neofytos Rodosthenous$^\ddagger$}
\thanks{$^\ddagger$Department of Mathematics, University College London, Gower St, London WC1E 6BT, UK. Email: n.rodosthenous@ucl.ac.uk}

\numberwithin{equation}{section}

\begin{document}

\begin{abstract} 
The socioeconomic impact of pollution naturally comes with uncertainty due to, e.g., current new technological developments in emissions' abatement or demographic changes. On top of that, the trend of the future costs of the environmental damage is unknown: 
Will global warming dominate or technological advancements prevail? 
The truth is that we do not know which scenario will be realised and the scientific debate is still open. 
This paper captures those two layers of uncertainty by developing a real-options-like model in which a decision maker aims at adopting a once-and-for-all costly reduction in the current emissions rate, when the stochastic dynamics of the socioeconomic costs of pollution are subject to Brownian shocks and the drift is an unobservable random variable. 
By keeping track of the actual evolution of the costs, the decision maker is able to learn the unknown drift and to form a posterior dynamic belief of its true value. 
The resulting decision maker's timing problem boils down to a truly two-dimensional optimal stopping problem which we address via probabilistic free-boundary methods and a state-space transformation. 
We completely characterise the solution by showing that the optimal timing for implementing the emissions reduction policy is the first time that the learning process has become ``decisive'' enough; that is, when it exceeds a time-dependent percentage. This is given in terms of an endogenously determined threshold function, which solves uniquely a nonlinear integral equation.
We numerically illustrate our results, discuss the implications of the optimal policy and also perform comparative statics to understand the role of the relevant model's parameters in the optimal policy.
\end{abstract}

\maketitle






 
\section{Introduction}
\label{sec:intro}

In 2006, the economist Nicholas Stern presented the famous report \textit{The Economics of Climate Change: The Stern Review} that was commissioned by the British government. 
It called for immediate and strong actions to reduce greenhouse gas emissions to prevent significant losses in global gross domestic product (GDP) with considerably cheaper actions 
(see \cite[Summary of Conclusions]{Stern}).
Even though many experts did not agree with all the assumptions and conclusions made in the report (see \cite{Nordhaus1}, \cite{Tol2} and \cite{Weitzman} for an overview), the {\it Stern Review} has contributed considerably to raise awareness of global warming.
Sustainability is nowadays one of the most important topics. 
State unions and governments are making bold policy statements, companies across a spectrum of industries are making their own policy moves, and this is only the beginning. The European Union's climate target plan is at least a 55\% reduction in greenhouse gas (GHG) emissions by 2030, the UK's is 80\% by 2050 (against 1990), and Germany's is 38\% by 2030 (against 2005), as set in the Effort Sharing Regulation (ESR). 
Industries (e.g.\ consumer packaged goods, etc.) will have to keep a check on emissions and meet the given standards, via their own policy moves towards net zero targets (e.g.\ Nestl\'e by 2050, Unilever by 2029, etc.).  
The public's principles are also aligned towards this direction, as ``today’s consumer asks even more than before for sustainability'' (Mark Schneider, CEO Nestl\'e) and ``sustainability is the top issue for investors" (Larry Fink, CEO BlackRock), which puts even more emphasis towards achieving these targets.
However, more than 15 years after the {\it Stern Review}, many important political questions remain largely unanswered and the debate on climate policies is still convoluted, especially due to the uncertainty and the irreversibility inherently related to those actions. By using the words of \cite{Allenetal}, ``Global efforts to mitigate climate change are guided by projections of future temperatures. But the eventual equilibrium global mean temperature associated with a given stabilization level of atmospheric GHG concentrations remains uncertain, complicating the setting of stabilization targets to avoid potentially dangerous levels of global warming''.

%
%
%
%

To begin with, the standard {\it cost-benefit analysis} used by businesses and policy makers for decision making is inappropriate for environmental policies,
primarily due to the presence of uncertainty in the evolution of the ecosystem and its resulting social and economic impacts (see \cite{LaRiviere-etal}, \cite{Pindyck2} for an overview), as well as the involvement of two important kinds of irreversibility. 

The socioeconomic uncertainty stems from the fact that the damages and costs of environmental pollution and GHG emissions (e.g.\ carbon dioxide (CO\textsubscript{2}) -- the main pollutant driving global warming) are barely foreseeable. 
This is due to e.g.\ the diverse and complex effects of an increase in the average global temperature, such as rising sea levels, increasing frequency and intensity of catastrophic events, storms, hurricanes, heat waves, as well as decreasing the cold stress, reducing energy demand for heating. Moreover, GHG emissions are fundamental for the energy system, food production, etc., and their sources are in fact every company and household.

In terms of irreversibilities, on one hand, environmental damage can be partially or even completely irreversible. 
Consider for example CO\textsubscript{2}, which stays in the atmosphere for hundreds of years and its atmospheric concentration reduces very slowly,  
or the potentially permanent damages caused by an increased average temperature. 
Clearly, these kinds of irreversibilities imply a sunk benefit that is associated with early policy adoption.
On the other hand, there are also always sunk costs associated with policy adoption. 
For example, the loss of employment, GDP reductions and  significant investments in abatement equipment by companies to avoid pollution; opportunity costs that bias in favour of waiting for new information and delaying the policy adoption. The effects of uncertainty and these types of irreversibility are therefore ambiguous.

Nowadays, the literature on optimal pollution management is huge, so that any attempt of a review would necessarily lead to a non exhaustive list of contributions. 
We therefore focus solely on the branch of works dealing with optimal timing decisions in environmental economics, which is where our contribution lies. 
In this regard, an early influential contribution is \cite{Pindyck1}, which also provides an overview of former studies. 
This work studies how uncertainty over future costs and benefits of reduced environmental degradation interact with the irreversibility of the sunk costs associated with an environmental regulation, and the sunk benefits of avoided environmental degradation. 
The ways in which various kinds of environmental and socioeconomic uncertainties can affect optimal policy design are then discussed in \cite{Pindyck2}. 
More recently, the optimal timing and size of pollution reduction in polluted areas is considered in \cite{Lappi} and the carbon emissions reduction is considered in \cite{HLG21} from the viewpoint of individual companies aiming for the minimisation of costs from carbon taxes and investment costs, while \cite{Murto} focuses on the maximisation of production gains against investment costs. 
Furthermore, a model for the optimal switching decision from a fossil-fuelled to an electric vehicle, from an individual's perspective, is developed in \cite{Falboetal}.
In terms of applications of Bayesian learning methods from a real-options perspective, 
a model for evaluating energy assets and potential investment projects under dynamic energy transition scenario uncertainty is developed in \cite{FloraTankov}. This leads to irreversible investment problems (entry/exit problems) under Bayesian uncertainty, which are then solved numerically and for which empirical analysis is provided. 
An investment in a renewable energy project, when decision makers are uncertain about the timing of a subsidy revision and therefore update their belief in a Bayesian fashion is considered in \cite{Dalbyetal}; a detailed numerical analysis provides insights about the role of policy uncertainty in the case of fixed feed-in tariffs. 
The interplay between the inspections performed by a regulator and noncompliance disclosure by a production firm are investigated in \cite{Kim}. The model leads to a dynamic game where the regulator chooses the timing of inspections and the company whether it should disclose a random occurrence of noncompliance.

In this paper, we wish to introduce and analyse a model that captures the issues of uncertainty and irreversibilities in the timing problem faced by governments, regulatory bodies or unions of states, for adopting environmental policies, inspired by \cite{Pindyck1} (see also \cite{Lappi}, \cite{Pindyck0} and the discussion in \cite{Pindyck2}), while also introducing additional uncertainty around the future social and economic costs of pollution, which may be largely unpredictable.
The main goal is to rigorously investigate how the considered increased economic uncertainty interacts with irreversibilities in the decision of when to optimally adopt the policy. 

In the course of this, we take the point of view of a social planner that faces a real-options-like irreversible investment decision with sunk cost $I$ for the once-and-for-all reduction in the current emissions with rate $E>0$ to a smaller rate $\hat E \geq 0$. 
This reflects the fact that major new environmental policies are unlikely to be revised often. 
We assume that the pollution stock (e.g.\ the average atmospheric concentration of CO\textsubscript{2}), modelled in the spirit of \cite{Nordhaus2} (see also \cite{Pindyck1}), generates damage that can be measured and put into monetary terms.
Therefore, there exists a stochastic process $X=(X_t)_{t\geq 0}$ modelling the random evolution of social and economic costs, associated to each unit of pollution stock $P=(P_t)_{t\geq 0}$ (or $\hat{P}=(\hat{P}_t)_{t\geq \tau}$ after the policy adoption at time $\tau$).
The social planner's aim is to choose an optimal (random) time $\tau$ to adopt the policy, so that the reduced future costs $(X_t \hat{P}_t)_{t \geq \tau}$ are closer to socially optimal levels, by incurring a sunk investment cost $I$ (e.g.\ pollution abatement equipment) associated to the adopted environmental policy (e.g.\ carbon tax). 
This sunk cost creates an incentive to wait for new information, contrary to the desire for policy adoption that would decrease future costs, leading to an interesting trade-off for the social planner who wishes to minimise the overall future expected costs.

As in existing literature, the {\it first layer of economic uncertainty} comes through the stochastic fluctuations in the dynamics of the socioeconomic impact of pollution $X$. 
The socioeconomic costs however could be increasing on average over time (e.g.\ African crop yields could be reduced by up to 50\% due to climate change) if global warming dominates, or decreasing (e.g.\ more efficient agriculture due to high-developed farms with the same soil and climatic prerequisite, access to high-quality seeds, pesticides) if the yield gap is closed due to technological advancements \cite{Tol1}. The truth is that we do not know which scenario will be realised, but being aware of this kind of uncertainty, we introduce a {\it second layer of economic uncertainty} in our model by assuming that the social planner has only partial information about $X$. 
This reflects an {\it uncertainty over the uncertainty} (see \cite{Pindyck2} and also \cite{Barnett} for a general dynamic equilibrium model under Knightian uncertainty).   
The main purpose of this work is therefore to provide a new framework to deal with this extensive uncertainty over future social and economic costs of pollution in the optimal timing of environmental decisions.

Given that each of these uncertainties aggravates in time, especially over long time horizons, we firstly assume that $X$ is a geometric Brownian motion, which considers that the exacerbations are of exponential type. 
The unpredictable nature of future costs of environmental pollution is modelled by an unpredictable drift $\mu$ (expected/average future impact) which is considered random and non-observable by the social planner in our novel modelling approach. 
To be more precise, we let $\mu$ be a discrete random variable that can take two values $\mu \in\lbrace -\alpha,\alpha\rbrace$ for some $\alpha>0$, a setup that represents the most crucial situation in a tractable way.  
If $\mu=\alpha$, then the cost per unit of pollution increases exponentially on average, incentivising a rather early policy adoption to reduce emissions. 
However, if $\mu=-\alpha$, then the cost per unit of pollution decreases exponentially on average, potentially due to inventions and technological advancements tackling future environmental pollution, incentivising the delay of policy adoption. 
This indeed reflects the contrary dynamic stemming from pessimistic projections of increasing future costs of pollution, based on the slow global progress so-far (e.g.\ CDP report \cite{CDP} states that change is not happening at the scale required), political challenges and lack of global cooperation to take aligned actions (e.g.\ according to the SBTi report \cite{SBTi}, 50\% of companies are off track to meet their climate targets), versus the optimism for successfully reducing the socioeconomic costs of pollution by major business players. These include Tesla, whose mission is to accelerate the world’s transition to sustainable energy via electric vehicles, innovations in energy storage and solar technology (e.g.\ Solar Roof, Powerwall), Breakthrough Energy Ventures, whose aim is to provide the world affordable, abundant clean energy via breakthrough technologies in energy, transportation, and agriculture, as well as Microsoft, Alphabet Inc., Amazon, Virgin Group, etc., with investments in carbon removal technologies and clean energy projects.
 
The approach that we use to analyse the problem involves the introduction of Bayesian learning via the a posteriori belief process $\Pi_t=\mathbb{P}(\mu=\alpha\vert \mathcal{F}_t^X)$. 
The idea is that given the social planner's partial information on $X$, their belief about the true drift $\mu$ is updated continuously as new information arrives via the real-time observation of the evolution of socioeconomic impact of pollution $X$, given by its natural filtration $\mathcal{F}_t^X$  
(this technique goes back to \cite{Shiryaev1} in a different context). 
Even though the two layers of uncertainty, the stochasticity of the dynamics of $X$ and its unobservable random drift, are initially independent, they become correlated via the learning process. 
The two layers of uncertainty essentially take the form of the belief process $\Pi$ about the true drift and the updated dynamics of $X$, whose drift is now dictated by the observable process $\Pi$, and they are driven by a common noise. 
There are some very recent studies with a similar mathematical background, such as 
an investment timing project \cite{Decamps1}, an optimal dividend problem \cite{DeAngelis2} and a Dynkin game \cite{DeAngelis1}, when the drift of an asset or firm's revenues is random and unknown to the manager, 
and an inventory management problem with unknown demand trend \cite{FFR}.
The reason for this growing interest in such models is that often the drift term of the underlying random process is unknown to decision makers, and estimating this parameter is a challenging task. As we have stressed before, this uncertain nature appears also in the evaluation of adopting environmental policies. However, to the best of our knowledge, the complete rigorous treatment of such a novel feature has never appeared before in the literature of optimal timing problems in environmental economics.

Decision makers with partial observations thus need to decide an optimal strategy, while simultaneously learning (updating their beliefs about) the unknown uncertainty via $\Pi$.
The resulting formulation under this framework of increased uncertainty leads to a three-dimensional optimal stopping problem with an underlying state space $(X,P,\Pi)$. 
We firstly show rigorously that the problem can be reduced to a two-dimensional optimal stopping problem.
A similar dimensionality reduction from a two- to a one-dimensional setting was conjectured in \cite{Pindyck1} for the full information version of this problem.  
Such one-dimensional optimal stopping problems (as in \cite{Pindyck1}) can often be solved analytically via the traditional \textit{guess-and-verify-approach}. 
This is non-feasible in the resulting two-dimensional setting in our paper though, since explicit solutions are typically not available due to the problem's associated PDE variational formulation. 
The resulting novel problem is indeed considerably harder to analyse than its standard full information version. 
The methodology employed to deal with the resulting genuine two-dimensional problem with coupled diffusive coordinates, includes a combination of probabilistic techniques and a state-space transformation to achieve enough regularity of the value function and the complete characterisation of the optimal stopping strategy (see also \cite{DeAngelis1}, \cite{Decamps1}, \cite{FFR}, and \cite{JP17}).

Our main result is the proof of a fine regularity of the problem's value function, and, more importantly, of the complete and practically implementable characterisation of the optimal stopping time: It is optimal to reduce emissions when the estimate $\Pi$ of the unknown cost trend becomes ``decisive'' enough, i.e.\ exceeds the boundary $c(Z) \in (0,1)$, where $Z$ is a deterministic process (like a ``time'' coordinate, whose speed is however determined by the size of socioeconomic costs' volatility). In particular, we show that the continuous curve $c$ uniquely solves a nonlinear integral equation, 
which is a considerable generalisation of the full information case  whose optimal policy adoption time is the hitting time of a constant barrier.
We would like to stress that, besides its theoretical importance, this is also a fundamental step in inferring application-driven conclusions.
As a matter of fact, at any time $t\geq 0$, the optimal policy adoption depends solely on the decision maker's degree of certainty $\Pi_t$ about increasing future costs, exceeding a deterministic $100 c(Z_t)\%$--confidence level (see Section \ref{sec:mainthm} for detailed results). 
An interesting outcome is the fact that $t \mapsto c(Z_t)$ turns out to be increasing, reflecting the fact that as time passes and more information is revealed about the a priori unknown socioeconomic cost trend, the decision maker becomes more reluctant to rush a policy adoption and is willing to wait until the certainty (based on the information gathering) for an increasing cost trend is higher. 
This provides a rigorous quantitative way to capture this important real-life phenomenon. 

The theoretical results are complemented by a numerical analysis aimed at determining how the different model parameters influence the optimal decision policy. Amongst various findings, we observe that the two layers of uncertainty induce different effects on the expected optimal timing of policy adoption. While an increase in the volatility $\sigma$ of the socioeconomic cost process $X$ induces an increase in the expected optimal time of pollution reduction -- in line with the classical ``value of waiting'' paradigm in real options -- increasing the average rate $\alpha$ of increase/decrease of future socioeconomic costs $X$, the variances of both the learning process $\Pi$ and of the unknown trend of $X$ increase, with the effect that decision makers become more proactive and act earlier on average.

The outline of the rest of the paper is as follows.
In Section \ref{frame} we introduce the decision maker's optimal timing problem, where the two layers of economic uncertainty interact with the irreversibility of the emissions reduction choice. 
Then, in Section \ref{section4}, we derive the equivalent three-dimensional Markovian formulation of the problem via filtering techniques.
The main result and its implications are presented in Section \ref{sec:mainthm}, together with comparative statics analysis of how certain model parameters affect the optimal strategy. 
Section \ref{sec:proof} then provides a constructive proof of the main theorem. 
This is distilled through a series of subsections and intermediate results. 
In particular, in Section \ref{section4n}, we prove that the three-dimensional problem can be reduced to a truly two-dimensional one. 
In Section \ref{sec:Solution}, we thus provide preliminary properties of the value function of the considered optimal stopping problem and of the boundary separating the action and waiting regions. 
The complete characterisation of the optimal policy adoption timing is then achieved in Section \ref{sec:Smooth}, via a state-space transformation developed in Section \ref{section56}. 
Section \ref{sectionNumAlgo} describes the algorithm to numerically illustrate the barrier which characterises the optimal stopping time.
Finally, we present our conclusions in Section \ref{section6}, where we also discuss problems with a similar structure that can be treated by techniques analogous to those employed in the present paper and we present ideas for future research directions. 
Appendix \ref{app:proofs} concludes this work by collecting the proofs of technical results.

\section{The optimal timing problem with uncertainty over uncertainty}
\label{frame}

Let $(\Omega,\mathcal{F},\mathbb{P}_\pi)$ be a complete probability space, rich enough to accommodate a one-dimensional Brownian motion $(B_t)_{t\geq 0}$ and a discrete random variable $\mu$ taking values $-\alpha$ and $\alpha$, with probability $1-\pi$ and $\pi$, respectively. 
Formally, the probability space $(\Omega,\mathcal{F},\mathbb{P}_\pi)$ is constructed by taking $\Omega:=C([0,+\infty),\mathbb{R})\times\lbrace -\alpha,\alpha\rbrace$, for some $\alpha>0$, and $\mathcal{F}$ as a $\sigma$-algebra satisfying the usual conditions. The probability measure $\mathbb{P}_\pi$ is then given by the product measure $\mathbb{P}_\pi:=\mathbb{W}\otimes\Lambda(\pi)$, where $\mathbb{W}$ denotes the standard Wiener measure and $\Lambda(\pi):=(1-\pi,\pi)$ is the discrete probability measure that assigns probability $1-\pi$ on $-\alpha$ and probability $\pi$ on $\alpha$, for some fixed $\pi\in(0,1)$. Then, the couple $\left((B_t)_{t\geq 0},\mu \right)$ is a canonical element in $\Omega$.

Let $P=(P_t)_{t\geq 0}$ be some pollutant stock, e.g.~the average atmospheric concentration of CO\textsubscript{2}, that evolves over time according to the ordinary differential equation (ODE)
\begin{equation}\label{P}
dP_t=(\beta E-\delta P_t)dt, 
\qquad P_0=p>0,
\end{equation}
where $E>0$ denotes the current level of emissions, $\beta>0$ is a scale parameter and $\delta>0$ is the dissipation rate of the pollutant. 
In particular, the ODE \eqref{P} admits a closed form solution  for all $\delta>0$, given by 
\begin{align}
P_t^p&=\frac{\beta E}{\delta}\left( 1-e^{-\delta t}\right)+pe^{-\delta t} , \qquad \text{for all } t \geq 0. \label{po1}
\end{align}

\begin{remark}
It is worth mentioning that, the dynamics of the pollution stock in \eqref{P} are in line with the one used by Nordhaus \cite{Nordhaus2} to evaluate greenhouse gas reducing policies, in the context of climate change. 
Note that, \cite{Nordhaus2} assumed that social costs come from higher temperatures driven by an increasing atmospheric concentration of greenhouse gases, while here we will allow social costs to be generated by the pollution stock $P$ directly.  
\end{remark} 

In this framework, we assume that the instantaneous cost or society's disutility from pollution does not only depend on the current level of pollution stock $P$, but also on the current level of social and economic costs  $X=(X_t)_{t\geq 0}$ generated by a unit of pollution. 
Given that there is uncertainty around $X$ and its real-life exacerbation in long time horizons, we model it as an It\^o-diffusion evolving according to the stochastic differential equation (SDE)
\begin{equation}\label{XX}
dX_t=\mu X_tdt+\sigma X_t dB_t, \qquad X_0=x>0,
\end{equation}
where $\mu \in \R$ denotes the average rate of increase/decrease of future costs, the process $(B_t)_{t\geq 0}$ models all the exogenous shocks affecting the environmental sustainability (e.g. related technological achievements and new scientific discoveries in related fields, or the lack of means to tackle global warming) and the volatility $\sigma>0$ denotes their extend. 

The main novelty of our model is that future social and economic costs of pollution are considered unpredictable, which plays a crucial role in the debate on environmental policies. 
Hence, we assume that the social planner has only partial information about the level $X$ of socioeconomic impact of pollution and cannot observe (estimate) the random rate $\mu$ of expected (average) future costs. 
This reflects the additional uncertainty around the instantaneous trend of technological advances and socioeconomic impacts of pollution (for the importance of uncertainty over existing uncertainty, see e.g.~\cite{Pindyck2}).
Moreover, the discrete random variable $\mu \in\lbrace -\alpha,\alpha\rbrace$, for some $\alpha>0$, represents in a tractable way the most crucial situation, which is in line with recent contrasting opinions of experts. 

The social planner has a prior belief on that $\mu=\alpha$, given by some fixed $\pi\in (0,1)$, and can only observe the evolution of the overall socioeconomic impact of pollution $X$. 
The process $X^x=(X_t^x)_{t\geq 0}$ satisfying the dynamics in \eqref{XX} is therefore a geometric Brownian motion whose drift depends on the unobservable random variable $\mu$ and it is such that
\begin{align}
\mathbb{E}^{\mathbb{P}_\pi} \left[ X_t^x\right]
&=\mathbb{E}^{\mathbb{P}_\pi}\Big[ xe^{\sigma B_t-\frac{\sigma^2}{2}t}\left(\eins_{\lbrace\mu=\alpha\rbrace}e^{\alpha t}+\eins_{\lbrace\mu=-\alpha\rbrace}e^{-\alpha t}\right)\Big]\nonumber\\
&=\mathbb{E}^{\mathbb{W}}\Big[ xe^{\sigma B_t-\frac{\sigma^2}{2}t}\Big]\,\mathbb{E}^{\Lambda(\pi)}\left[ \left(\eins_{\lbrace \mu=\alpha\rbrace}e^{\alpha t}+\eins_{\lbrace\mu=-\alpha\rbrace}e^{-\alpha t}\right)\right]
=x\left(\pi e^{\alpha t}+(1-\pi)e^{-\alpha t}\right),\label{expectation}
\end{align}
where $\mathbb{E}^{\mathbb{P}_\pi}\left[\,\cdot\,\right], \mathbb{E}^{\mathbb{W}}\left[\,\cdot\,\right], \mathbb{E}^{\Lambda_\pi}\left[\,\cdot\,\right]$ denote the expectations under the probability measures $\mathbb{P}_\pi, \mathbb{W}$, $\Lambda(\pi)$, respectively. 
Notice that the second equality follows due to the independence of the process $(B_t)_{t\geq 0}$ and the random variable $\mu$.

\begin{remark}[Full information] \label{rem:full}
The case of a known constant rate $\mu$ (in practice, estimateable), such that the economic uncertainty is fully observable and derived solely from the diffusion term has been considered in \cite{Pindyck1}.
In this case, the process $X$ defined by \eqref{XX} is a geometric Brownian motion and its closed form solution is given by 
\begin{equation*}
X_t^x=x\exp\Big\{\Big(\mu-\frac{1}{2}\sigma^2\Big)t+\sigma B_t\Big\} , \qquad \text{for all } t \geq 0. 
\end{equation*}
\end{remark}

Overall, if at time $t\geq 0$, the level of pollution is $p$ and the marginal social and economic cost is $x$, then the cost generated by the environmental pollutant is $x p$.
Taking this into account, we consider a social planner whose target is to choose a (random) time $\tau$ at which an environmental policy should be adopted in order to reduce the emissions rate from $E$ to some lower $\hat E$. 
In this paper, we consider $\hat E=0$ without loss of generality. 
Hence, the pollutant stock after the policy adoption, denoted by $(\hat{P}_t)_{t\geq \tau}$, will follow the dynamics 
\begin{equation}\label{p}
d\hat{P}_t=(\beta \hat{E} - \delta\hat{P}_t ) dt = -\delta \hat{P}_t dt , \qquad \text{for all } t > \tau, \qquad  \hat{P}_\tau=P_\tau.
\end{equation}
In particular, the ODE \eqref{p} admits a closed form solution given by 
\begin{align}
\hat{P}^p_t&= P^p_\tau \, e^{-\delta (t-\tau)} , \qquad \text{for all } t\geq\tau.\label{po2}
\end{align}
Finally, it is natural to assume that any environmental policy adoption yields other societal and economic costs, e.g. due to 
loss of employment, reductions in the GDP, costly investments in abatement equipment. 
We assume that this investment cost is completely sunk and given by the constant $I>0$. 

Given a constant discount rate $r>0$ and any initial values $x,p>0$, the social planner's objective is to find a stopping time $\tau$ of the filtration ${\mathcal{F}}_t^X$ generated by $X$ (representing the information flow generated by observing the actual evolution of the socioeconomic costs of pollution), at which it is optimal to spend the investment costs $I$ in order to permanently reduce the emissions from rate $E$ to $\hat E$. 
This target can be formulated via a (non-Markovian) optimal stopping problem over an infinite time horizon given by 
\begin{equation} \label{valf1}
\inf_{\tau\geq 0}\,\mathbb{E}^{\mathbb{P}_\pi}\left[ \int_0^\tau e^{-rt} X_t^x P_t^p dt + e^{-r\tau}I + \int_\tau^\infty e^{-rt} X_t^x \hat{P}^p_t dt\right],
\end{equation}
where the infimum is taken over all stopping times $\tau$ of the process $(X_t^x)_{t\geq 0}$.
Notice that, the first integral in the expectation in \eqref{valf1} represents the cumulative costs until the policy is adopted, while the second one the cumulative costs after the policy adoption. 

\begin{remark} \label{rem:PdB}
The mathematical analysis in this paper applies also in the case of a stochastically evolving stock of pollutants, namely when
\begin{equation*} 
dP_t=(\beta E - \delta P_t)dt + \eta d{\widetilde B}_t, 
\qquad P_0=p>0,
\end{equation*}
where the parameters $E, \beta, \delta$ are as in \eqref{P}, while $(\widetilde B_t)_{t\geq 0}$ is a Brownian motion (independent of $B$) modelling the shocks affecting the atmospheric stock of pollutants and the volatility $\eta>0$ denotes their extend (see \cite{AthanassoglouXepapadeas}, Section 3 in \cite{Pindyck0} and Section 5 in \cite{Pindyck1}, among others).

Such dynamics will neither interfere with the learning process of the decision maker nor affect the analysis resulting to the two-dimensional problem \eqref{valf2} (cf.\ Section \ref{section4}). 
The only difference is that the expected values of the stock of pollutants should be used in the calculations, instead of the explicit expressions \eqref{po1} and \eqref{po2}. 
All subsequent analysis of the resulting two-dimensional problem should be identical (see also Example 1 in Section \ref{Sec:similar}). 
\end{remark}

\subsection{Markovian formulation of Problem \eqref{valf1}}
\label{section4}

In order to solve the problem \eqref{valf1}, we first observe that a Markovian reformulation is needed, since the unpredictable and non-observable nature of $\mu$ implies that $X$ is not a Markov process. The first step of our forthcoming analysis is thus to express the optimal stopping problem \eqref{valf1} in a Markovian framework.

The social planner's information is modelled by the filtration $\mathcal{F}^X=(\mathcal{F}_t^X)_{t\geq 0}$ generated by $X$ and augmented with $\mathbb{P}_\pi$-null sets, which is right-continuous (cf. \cite[Theorem 2.35]{Bain}). 
The social planner can then update their belief on the true value of $\mu$ according to new information as it gets revealed. 
In other words, by relying on filtering techniques (see, e.g. \cite{Liptser}), we can define the social planner's Bayesian learning 
process $\Pi = (\Pi_t)_{t\geq 0}$ on $(\Omega,\mathcal{F},\mathbb{P}_\pi)$ as the ${\mathcal{F}}^X$-càdlàg martingale  
\begin{equation}
\label{eq:defPi}
\Pi_t:=\mathbb{P}_\pi\big(\mu=\alpha \,\big|\,{\mathcal{F}}_t^X\big), \qquad \Pi_0=\pi \in (0,1) .
\end{equation} 
Then, by \cite[Section 4.2]{Shi}, the process $(X^{x,\pi},\Pi^\pi)$ uniquely solves the SDE 
\begin{align} \label{XPi}
\begin{cases}
 dX_t^{x,\pi}=\left(\alpha \Pi_t^\pi-\alpha(1-\Pi_t^\pi)\right)X_t^xdt+\sigma X_t^x d{W}_t, &\quad X_0^{x,\pi} = x > 0,\\ 
d\Pi_t^\pi=\frac{2\alpha}{{\sigma}}\Pi_t^\pi(1-\Pi_t^\pi)d{W}_t, &\quad \Pi_0^\pi = \pi \in (0,1),
\end{cases}
\end{align}
where $W=(W_t)_{t\geq 0}$ is the so-called innovation process, which is defined by 
\begin{equation*}
W_t := \int_0^t \frac{\mu+\alpha-2\alpha\Pi^{\pi}_s}{\sigma} \, ds + B_t,
\end{equation*}
and it is an ${\mathcal{F}}^X$-adapted Brownian motion under the probability measure $\mathbb{P}_\pi$. 
We note that the process $(X^{x,\pi},\Pi^\pi)$ is also ${\mathcal{F}}^{{W}}$-adapted and, thus, we have ${\mathcal{F}}^X={\mathcal{F}}^{{W}}$. 
It is clear that $(X^{x,\pi}, \Pi^{\pi})$ is now a Markov process under this new formulation. 

Consider now a new probability space $(\widetilde{\Omega}, \widetilde{\mathcal{F}}, \widetilde{\mathbb{P}})$, on which we define a Brownian motion $\widetilde{W}$, adapted to its natural filtration ${\mathcal{F}}^{\widetilde{W}}$, augmented by $\widetilde{\mathbb{P}}$-null sets of $\widetilde{\mathcal{F}}$. 
On such a space, let $(\widetilde{X},\widetilde{\Pi})$ evolve as in \eqref{XPi}, but with $W$ replaced by $\widetilde{W}$. Since \eqref{XPi} admits a unique strong solution, then
$$\text{Law}_{\mathbb{P}_{\pi}}(X, \Pi, \tau) = \text{Law}_{\widetilde{\mathbb{P}}}(\widetilde{X}, \widetilde{\Pi}, \widetilde{\tau}),$$
where $\widetilde{\tau}$ is an ${\mathcal{F}}^{\widetilde{W}}$-stopping time.
Then, by means of the tower property in the expectation 
of \eqref{valf1}, and using \eqref{eq:defPi}--\eqref{XPi} and the aforementioned equality in law, the optimal stopping problem \eqref{valf1} becomes
\begin{equation}
\label{eq:OStilde}
\inf_{\widetilde{\tau}\geq 0}\,\widetilde{\mathbb{E}}\Big[\int_0^{\widetilde{\tau}} e^{-rt} \widetilde{X}_t^{x,\pi} {P}_t^{p} dt + e^{-r\widetilde{\tau}}I + \int_{\widetilde{\tau}}^\infty e^{-rt} \widetilde{X}_t^{x,\pi} \hat{P}_t^{p} dt\Big].
\end{equation}
From now on, with a slight abuse of notation, we shall write $(X,\Pi,\tau)$ as well as $(\Omega, \mathcal{F}, \mathbb{P}, \mathbb{E})$ instead of $(\widetilde{X},\widetilde{\Pi}, \widetilde{\tau})$ and $(\widetilde{\Omega}, \widetilde{\mathcal{F}}, \widetilde{\mathbb{P}}, \widetilde{\mathbb{E}})$, respectively.

Then, with regards to the Markovian nature of \eqref{eq:OStilde}, given any initial values $(x,p,\pi)\in\mathbb{R}_+\times\mathbb{R}_+\times (0,1)$, we consider the optimal stopping problem 
\begin{equation}\label{partprob1}
V(x,p,\pi):=\inf_{\tau\geq 0}\,\mathbb{E}_{x,p,\pi}\Big[\int_0^\tau e^{-rt} {X}_t P_t dt + e^{-r\tau}I + \int_\tau^\infty e^{-rt} {X}_t \hat{P}_t dt\Big],
\end{equation}
where $\mathbb{E}_{x,p,\pi} [\,\cdot\,] :=\mathbb{E}[\,\cdot\,\big| \, X_0=x, P_0=p, \Pi_0=\pi]$ and the infimum is taken over all ${\mathcal{F}}^{X}$-stopping times $\tau$. 
Solving \eqref{partprob1} then consists of finding the optimal timing $\tau^*$ for adoption the environmental policy that achieves the minimum overall expected socioeconomic costs $V$.

\section{The main result and its implications}
\label{sec:mainthm}

Our main result provides a complete characterisation of the optimal policy adoption time and it is provided in the theorem below. 
In order to consider the potential optimality of the immediate emissions reduction policy or its perpetual postponement, which are certainly plausible choices in environmental economics (especially in environmental policy adoption discussions), we make the only standard assumption that $r>\alpha$.
The rest of this paper is then devoted to develop a constructive proof of such a theorem.

\begin{theorem} 
\label{thm:mainthm}
Assume $r>\alpha$, recall the Bayesian learning process $\Pi^{\pi}$ defined by \eqref{eq:defPi} and denote its transition density by $p_t(\pi,\pi')$, for $(\pi,\pi')\in (0,1)^2$, define the auxiliary ``time-coordinate'' process 
\begin{align*}
Z^z_t &:=z + \frac{1}{2}\sigma^2 t, \quad t\geq 0, \quad  \text{with} \quad z=z(x,\pi):= \frac{\sigma^2}{2\alpha}\ln(\frac{\pi}{1-\pi}) - \ln(x),
\end{align*}
and introduce the continuous, nondecreasing function
$$
m(z):=\inf\{\pi\in (0,1)\,|\, q(z,\pi) <0\}, \quad z \in \mathbb{R},
$$
where 
\begin{align*} 
&q(z,\pi) := \beta E e^{-z} \Bigl(\frac{\pi}{1-\pi}\Bigr)^\frac{\sigma^2}{2\alpha} \Big( \big( (\alpha-r)\theta+2\alpha\rho \big) \pi - (\alpha+r)\rho \Big)+rI,
\quad (z,\pi) \in \R  \times (0,1), \\
&\theta :=\frac{2\alpha(2r+\delta)}{(r-\alpha)(r+\alpha)(r+\delta-\alpha)(r+\delta+\alpha)}  
\quad \text{and} \quad 
\rho :=\frac{1}{(r+\delta+\alpha)(r+\alpha)}.
\end{align*}

Then, aiming at the minimisation of the overall socioeconomic costs of pollution in \eqref{partprob1} for any $(z,\pi) \in \mathbb{R}\times (0,1)$, it is optimal to adopt the emissions reduction policy at the stopping time
$$\tau^{*}=\tau^*(z,\pi):=\inf\{t\geq 0:\, \Pi^\pi_t \geq c(Z^z_t)\},$$
where $c:\mathbb{R} \to [0,1]$ is the unique continuous nondecreasing solution to the integral equation
\begin{align} \label{solc}
0 = \int_0^{\infty} e^{-rt} \left( \int_{0}^{c(z+\frac12 \sigma^2 t)} \hspace{-3mm}q\big(z+\tfrac12 \sigma^2 t, \pi'\big) p_t(c(z),\pi') d\pi' \right) dt ,
\end{align}
such that $c(z) \geq m(z)$, for all $z\in \mathbb{R}$.
\end{theorem}

Besides the theoretical interest of proving Theorem \ref{thm:mainthm} and completing the analysis of our problem, this result provides a way to numerically implement our theoretical findings and therefore to understand the role of various model parameters in the optimal strategy.

We firstly notice that, the decision maker learns about the true value of the expected future cost trend $\mu$ via the learning process $\Pi_t \in [0,1]$, for $t \geq 0$. The latter process begins at time $0$ from the decision maker's initial belief $\pi \in (0,1)$ about $\mu=\alpha$, and it is then updated continuously as new information arrives via the real-time observation of the evolution of socioeconomic impact of pollution $X$, given by its natural filtration $\mathcal{F}_t^X$. 
The emissions reduction policy should be optimally adopted at $\tau^*$, namely as soon as the learning process $\Pi_{\tau^*}$ exceeds a deterministic threshold $c(Z_{\tau^*}) \in [0,1]$. 
Essentially, the latter can be viewed as a $100 c(Z_t)\%$--confidence level, for $t \geq 0$, which triggers the policy adoption when exceeded by the decision maker's degree of certainty $\Pi_t$ about increasing future costs $\mu=\alpha$, at the optimal stopping time $t = \tau^*$.
In other words, the decision maker should optimally adopt the emissions reduction policy as soon as they are ``confident enough'' that the socioeconomic costs of pollution are on average increasing.
 
One of the powers of the aforementioned result is that the optimal policy is completely characterised in a way that does not involve explicitly the stock of pollutants process $P$ and its socioeconomic costs $X$, but in a way that only their driving parameters are involved. 
To be more precise, even though the desired confidence level $100 c(Z^z_t)\%$ is deterministic, since it is driven by the auxiliary ``time-coordinate'' process $Z^z_t = z + \frac{1}{2}\sigma^2 t$, it is important to note that the speed with which it increases depends on $\sigma$, namely how volatile the socioeconomic costs of pollution are, cf.\ \eqref{XPi}. 
At the same time, the confidence level's starting value 
$$100 c(z)\% = 100 \, c \Big(\frac{\sigma^2}{2\alpha}\ln \Big(\frac{\pi}{1-\pi} \Big) - \ln(x) \Big)\%$$ 
depends on all parameters of $X$ and the decision maker's a priori belief $\pi$ about a high future impact $\mu=\alpha$ of costs. 
In addition, we observe that all these parameters together with the ones driving the evolution of $P$ are part of the integral equation whose solution is the actual function $c$, hence the shape of the confidence level function $z \mapsto 100 c(z)\%$ depends on all model parameters -- this is the case also for the lower bound $m(z)$ of admissible confidence levels $100 c(z)\%$, which is defined via the critical quantity $q(z,\pi)$.
Further analysis of the optimal policy result in Theorem \ref{thm:mainthm}, the implications of the monotonicity of $c$ and the effect of parameters on such a policy are provided in the subsequent discussion.

\subsection{Decision making phases}
\label{sec:phases}

\begin{figure} 
\includegraphics[width=10cm]{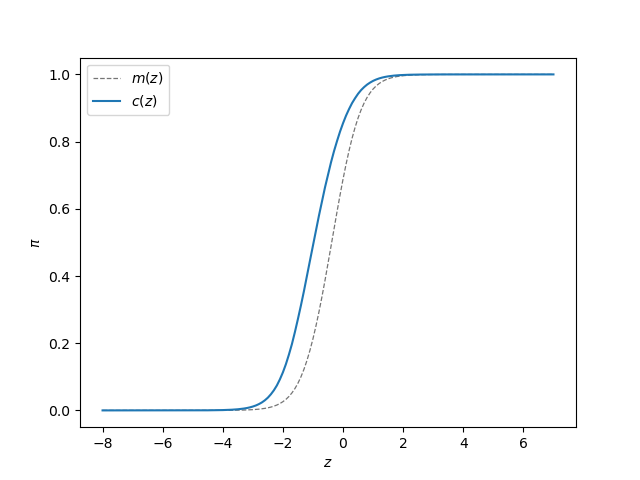}
\caption{A numerical calculation of the boundary function $z \mapsto c(z)$ solving \eqref{solc} and dominating the threshold $z \mapsto m(z)$ defined in Theorem \ref{thm:mainthm}.}
\label{c>m}
\end{figure}

In Figure \ref{c>m},  we plot the boundary function $z \mapsto c(z)$ and the lower threshold $z \mapsto m(z)$, as defined in Theorem \ref{thm:mainthm}. The parameters used, as well as the description of the numerical algorithm, are collected in Section \ref{sectionNumAlgo}.

The first main conclusion one can draw from Figure \ref{c>m} is that, based on the configuration of model parameters defining the initial value $z=z(x,\pi):= \frac{\sigma^2}{2\alpha}\ln(\frac{\pi}{1-\pi}) - \ln(x)$ of the process $Z$, introduced in Theorem \ref{thm:mainthm}, we may start our observations of the state space process $(Z,\Pi)$ from either one of three types of regions (e.g. for $z<-3.5$, $-3.5<z<2$, $z>2$ in Fig.~\ref{c>m}). 

In particular, the optimal timing problem for adopting an emissions reduction policy can experience the following three {\it Phases}:
\begin{itemize}
\item[(I).] \it Immediate adoption.
If $z$ is ``relatively small''  -- that is, if the initial socioeconomic cost of pollution $x$ is relatively large with respect to the adjusted likelihood ratio $(\frac{\pi}{1-\pi})^{{\sigma^2}/({2\alpha})}$ -- the initial belief $\pi$ of an increasing trend of future costs would be above $c(z)$, which would most likely require the optimal {immediate adoption} of the policy -- unless we are absolutely certain of a decreasing future cost trend, i.e.\ $\pi \approx 0$. 

\item[(II).] \it Dynamic decision making.
If $z$ takes a ``relatively intermediate'' value, we are in a {dynamic decision making phase}, where we decide to adopt the policy -- while learning the unknown future cost trend -- when the stochastically-evolving learning process $\Pi^{\pi}_\tau$ exceeds the critical deterministically-evolving threshold $c(Z^z_\tau)$ at some time $\tau \geq 0$.

\item[(III).] \it Never adopt.
If $z$ is ``relatively large'' -- that is, if the initial socioeconomic cost of pollution $x$ is relatively small compared to the adjusted likelihood ratio $(\frac{\pi}{1-\pi})^{{\sigma^2}/({2\alpha})}$ -- or if we start from Phase (II) and $\Pi^{\pi}_t$ remains below the increasing (in time) threshold $c(Z^z_t)$, we end up in this third phase, where we most likely {never adopt} the policy -- unless we are absolutely certain of an increasing future cost trend, i.e.\ $\Pi^{\pi}_\tau \approx 1$ at some time $\tau \geq 0$. 
\end{itemize}

In order to explore further the most interesting Phase (II) appearing in Figure \ref{c>m}, we can make two more observations:
(a) While learning in Phase (II), i.e.\ as time passes and $Z_\cdot$ increases, implying that $c(Z_\cdot)$ increases (see also properties of $c$ in Theorem \ref{thm:mainthm}), the decision maker requires a higher certainty about increasing future costs, in order to optimally adopt the policy; 
(b) The duration of time for which the observation process $(Z,\Pi)$ could stay in Phase (II), before adopting the policy or moving to Phase (III) and becoming ``too late'' for adopting the policy, is driven strongly by the value of the socioeconomic costs' volatility $\sigma$, since it determines the speed of the time-coordinate process $Z$ (see Theorem \ref{thm:mainthm}).

\subsection{Sensitivity analysis: Optimal timing of environmental policy adoption}
\label{SAtau}

Next, we aim at obtaining further results in terms of the sensitivity of the critical belief threshold $c(z)$ and the expected optimal timing for adopting the environmental policy $\E[\tau^* | \tau^* < \infty]$ with respect to several important model parameters.

\begin{figure} 
\includegraphics[width=8.1cm]{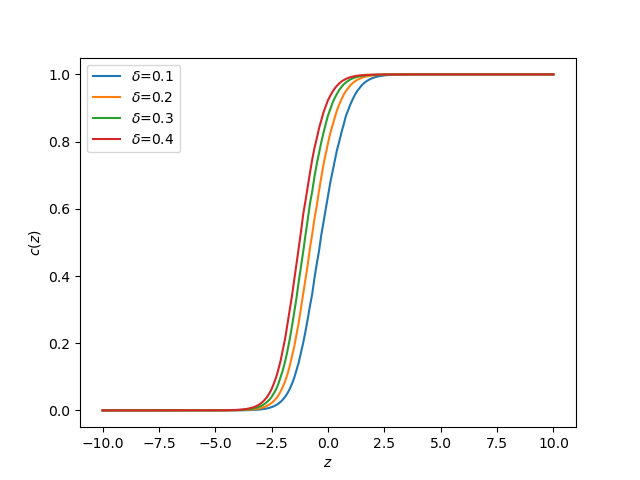}
\includegraphics[width=8.1cm]{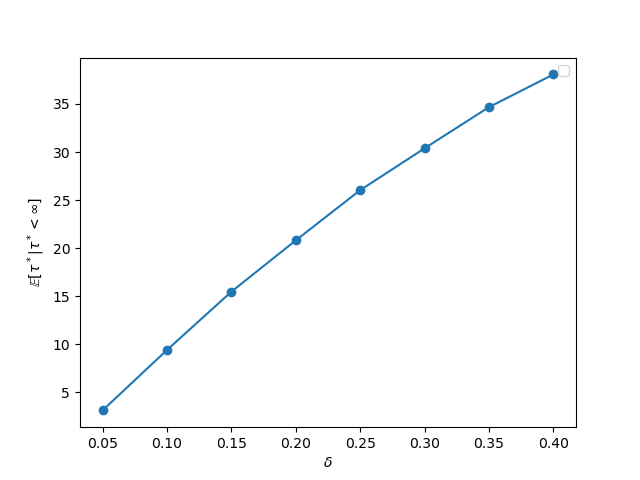}
\caption{A numerical calculation of the boundary function $z \mapsto c(z)$ solving \eqref{solc} with respect to different dissipation rates $\delta$ of the pollutant stock from the atmosphere, and the expected time to the policy adoption $\E[\tau^* | \tau^* < \infty]$ as a function of $\delta$.}
\label{figd}
\end{figure}

We see from Figure \ref{figd} that the critical belief threshold $c(z)$ is increasing with higher dissipation rates $\delta$ of the pollutant stock from the atmosphere. 
That is, decision makers become more reluctant to adopt the policy if the pollutant emitted will dissipate at a faster rate, irrespective of their actions, which also results in the delay of the optimal policy adoption observed in the right-hand panel. 
Essentially, they require their belief in an increasing future cost trend to reach a higher certainty level in order for the adoption of an emission reduction policy to be optimal.  
However, if the dissipation rate of pollutants is slow, then the importance of their actions increases, making the policy adoption optimal even for lower confidence levels about the future evolution of socioeconomic costs.

\begin{figure} 
\includegraphics[width=8.1cm]{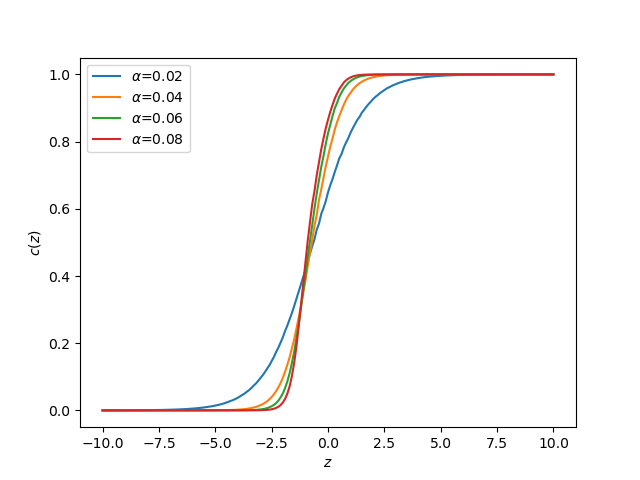}
\includegraphics[width=8.1cm]{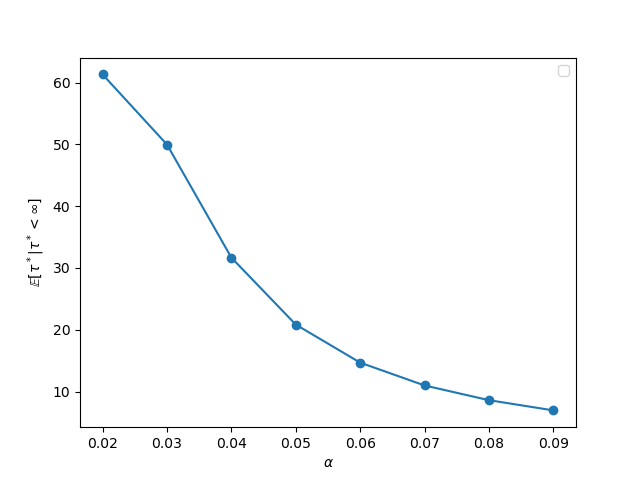}
\caption{A numerical calculation of the boundary function $z \mapsto c(z)$ solving \eqref{solc} with respect to different sizes $\alpha$ (resp. $-\alpha$) of the average rate of increase (resp.\ decrease) of future socioeconomic costs of pollution, and the expected time to the policy adoption $\E[\tau^* | \tau^* < \infty]$ as a function of $\alpha$.}
\label{figa}
\end{figure}

In Figure \ref{figa}, we observe an interesting phenomenon. 
There is no clear monotonicity of the critical belief threshold $c(z)$ with respect to changes in the absolute value $\alpha$ of the average rate of increase/decrease of future socioeconomic costs of pollution.
It seems though that, if we have more extreme alternative scenarios for the trend of future socioeconomic costs, i.e.\ higher $\alpha$-values, then the dynamic decision making phase (see Phase (II) in Section \ref{sec:phases}) shrinks in terms of time, and makes it more likely for the decision maker to end up in either an immediate policy adoption or never adopting the policy (both Phase (I) and Phase (III) in Section \ref{sec:phases}). 
Essentially, this provides a quantitative way to capture the fact that, as  the alternatives diverge from each other, the decision maker can learn sooner, i.e.\ after a relatively shorter time period, whether it is optimal to adopt the policy or not. 
Contrary, when $\alpha$ decreases and the two alternatives come closer to each other, the dynamic decision making phase is extended, so the decision maker requires a larger time-window to learn the unknown costs, while examining whether a policy adoption would be optimal.

Besides the aforementioned effects of a changing absolute value $\alpha$ of the average rate of increase/decrease of future socioeconomic costs of pollution on the decision maker's learning process, we can further conclude from the right-hand panel of Figure \ref{figa}, that the optimal timing of adopting the environmental policy decreases on average as $\alpha$ increases. 
In view of the dynamics of $\Pi^\pi$ in \eqref{XPi}, such higher values of $\alpha$ (also higher spread in alternative scenarios) result in the increase of the volatility in the decision maker's learning process about the unknown cost trend. 
This consequent increase in the {\it second layer of uncertainty} in our model (in the decision maker's prediction mechanism), results in more proactive decision makers, who are willing to bring forward their actions. 
Interestingly, this phenomenon is contrary to the philosophy of the ``value of waiting'', according to which decision makers usually become less proactive and are willing to postpone their actions in times of increased  uncertainty (see \cite{DiPi} and \cite{McD}, among others, for works presenting the classical ``value of waiting'' effect, and Section 4 in \cite{BT} for an instance where the ``value of waiting'' paradigm does not hold in the context of a costly sequential experimentation and project valuation). 
A result agreeing with the ``value of waiting'' is obtained in this model only in terms of the (more classical) {\it first layer of uncertainty}, which is presented below.  

\begin{figure} 
\includegraphics[width=8.1cm]{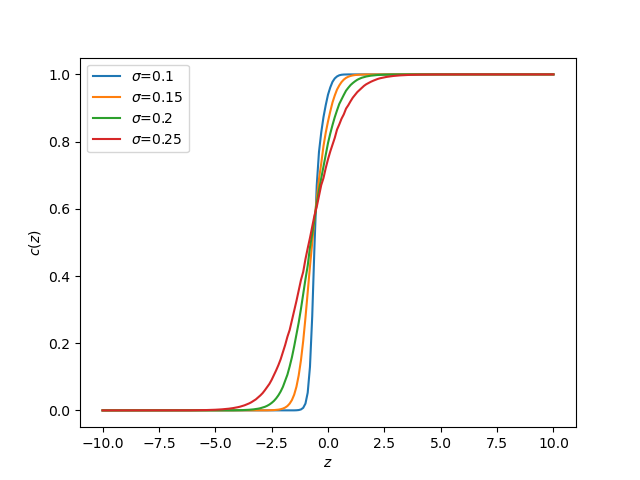}
\includegraphics[width=8.1cm]{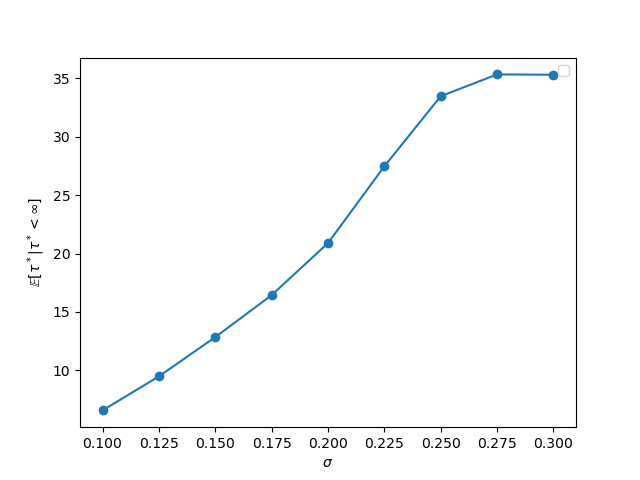}
\caption{A numerical calculation of the boundary function $z \mapsto c(z)$ solving \eqref{solc} with respect to different extents of volatility $\sigma$ in the socioeconomic costs of pollution, and the expected time to the policy adoption $\E[\tau^* | \tau^* < \infty]$ as a function of $\sigma$.}
\label{figs}
\end{figure}

We also observe in Figure \ref{figs} that the critical belief threshold $c(z)$ is clearly non-monotonic with respect to the volatility $\sigma$ of the socioeconomic costs of pollution. 
It seems that as the economic uncertainty around socioeconomic costs increases, i.e.\ higher $\sigma$-values, we get a larger dynamic decision making phase (see Phase (II) in Section \ref{sec:phases}). 
Therefore the decision maker requires a larger time-window to learn the unknown costs towards an optimal policy adoption.
On the contrary, in times of low uncertainty the dynamic decision making phase shrinks in terms of time, and makes it more likely for the decision maker to end up in either an immediate policy adoption or never adopting the policy (both Phase (I) and Phase (III) in Section \ref{sec:phases}). 
That is, with low volatility in socioeconomic costs, the future costs become less unpredictable and the current knowledge of decision makers would suffice to make an optimal decision to immediately or never adopt the policy. 
We can further conclude from the right-hand panel of Figure \ref{figs} that the decision makers become more reluctant to adopt the emissions reduction policy (optimal timing increases on average), if there is a higher {\it first layer of uncertainty} about the socioeconomic costs, which is consistent with the philosophy of the aforementioned ``value of waiting''.

\subsection{Sensitivity analysis: Optimal socioeconomic costs} 
\label{SAcost}

In what follows, we focus on the sensitivity of the expected socioeconomic costs, under the optimal policy adoption strategy of Theorem \ref{thm:mainthm}, with respect to several important model parameter.

\begin{figure} 
\includegraphics[width=8.1cm]{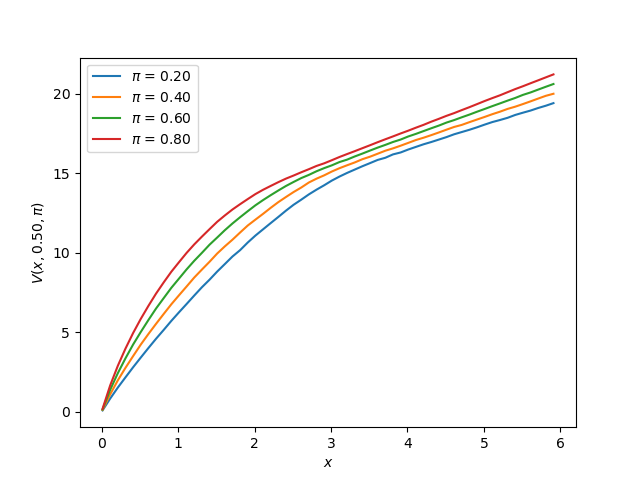}
\includegraphics[width=8.1cm]{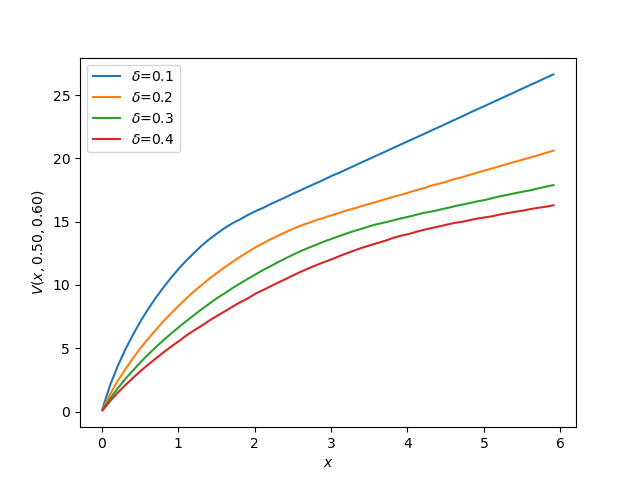}
\includegraphics[width=8.1cm]{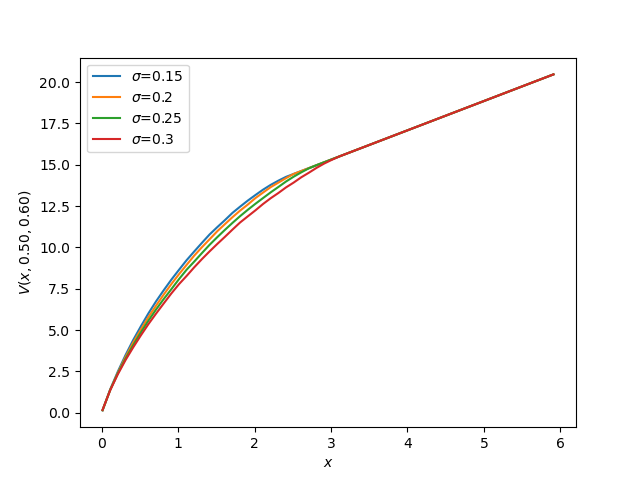}
\includegraphics[width=8.1cm]{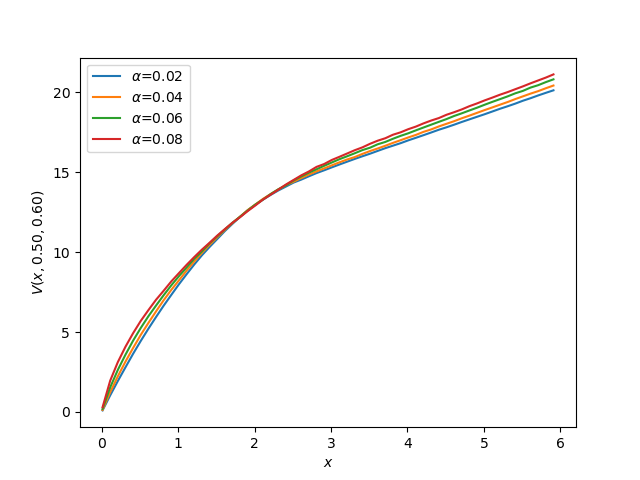}
\caption{A numerical calculation of the optimal expected socioeconomic costs $x \mapsto V(x,p,\pi)$ solving \eqref{partprob1} with respect to different: 
beliefs $\pi$ about an increasing average rate of future socioeconomic costs of pollution, 
dissipation rates $\delta$ of the pollutant stock from the atmosphere, 
extents of volatility $\sigma$ in the socioeconomic costs of pollution,
and sizes $\alpha$ of the average rate of increase/decrease of future socioeconomic costs of pollution.}
\label{figV}
\end{figure}

We begin with the more intuitive observations from the top panels of Figure \ref{figV}, namely, that the optimal expected socioeconomic costs are higher when the decision makers degree of certainty $\pi$ about increasing future socioeconomic costs of pollution increases, and lower when the dissipation rate $\delta$ of the pollutant stock from the atmosphere increases. 
In other words, the more certainty there is around increasing future costs of pollution and the slower the pollutant stock dissipates from the atmosphere, the higher the overall costs of pollution will be on average, even under the optimal environmental policy adoption.

Then, we observe some interesting features in the bottom panels of Figure \ref{figV}, in terms of the two layers of uncertainty in our model. 
Namely, the {\it first layer of uncertainty} in the socioeconomic costs $X$ per unit of pollution, which is increasing in $\sigma$, and the {\it second layer of uncertainty} due to the unknown trend of future costs, which is captured by the learning process $\Pi$ and is increasing in $\alpha$. 
In particular, the optimal expected socioeconomic costs are decreasing in $\sigma$ and increasing in $\alpha$. 
This implies that, when implementing the optimal policy adoption strategy, a higher {\it first layer of uncertainty} in our model works in the decision maker's favour (as it is usually observed in decision theory, e.g. financial option values are increasing in volatility), while, on the contrary, a higher {\it second layer of uncertainty} in our model works against the decision maker.

\subsection{Sensitivity analysis: Expected stock of pollutants until policy adoption}
\label{SApol}

\begin{figure} 
\includegraphics[width=8.1cm]{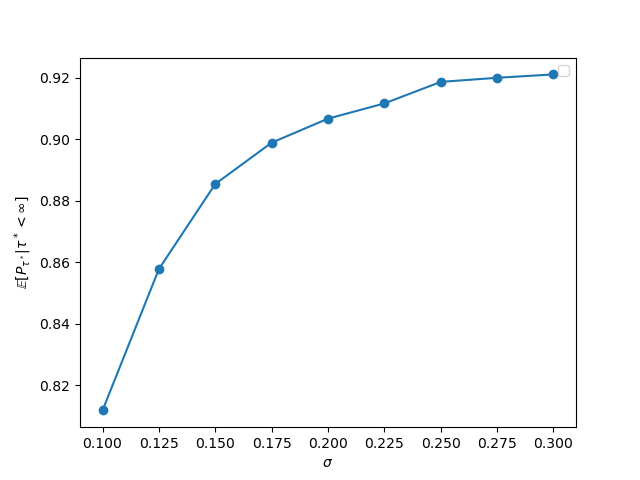}
\includegraphics[width=8.1cm]{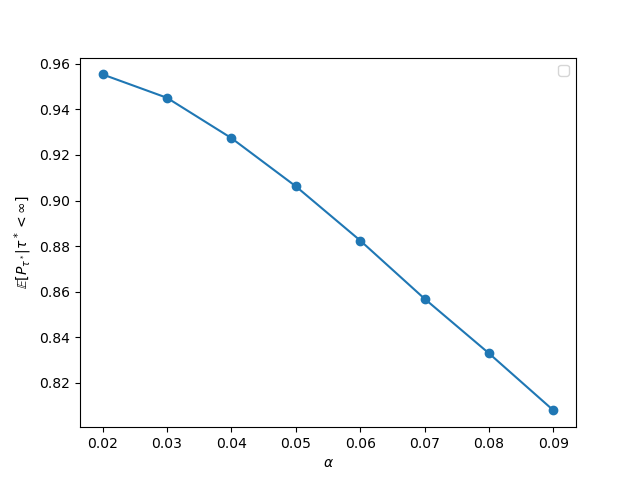}
\caption{A numerical calculation of the expected stock of pollutants in the atmosphere $\E[P_{\tau^*} | \tau^* < \infty]$ until the optimal adoption of the emissions reduction policy at time $\tau^*$, with respect to different extents of volatility $\sigma$ in the socioeconomic costs of pollution and sizes $\alpha$ of the average rate of increase/decrease of future socioeconomic costs of pollution.}
\label{figP}
\end{figure}

The interesting outcomes observed previously in Section \ref{SAcost}, in terms of the optimal cost of pollution against the two layers of uncertainty in our model, are extended in this section to their effect on the expected pollution until the optimal adoption of the emissions reduction policy at time $\tau^*$. 

To be more precise, we observe that, even though the overall socioeconomic cost of pollution is decreasing with the {\it first layer of uncertainty} in the socioeconomic costs $X$ per unit of pollution (see Figure \ref{figV} for $\sigma$ and Section \ref{SAcost}), the average stock of pollutants in the atmosphere up to the policy adoption is increasing in Figure \ref{figP}. 
Intuitively, this occurs due to the delay in the optimal policy adoption, observed in Figure \ref{figs} (see Section \ref{SAtau} for details), which results in giving more room for pollution stock to increase on average. 
A similar contradictory effect is also observed with respect to the {\it second layer of uncertainty} due to the unknown trend of future costs and the associated learning process $\Pi$, since the overall socioeconomic cost of pollution is increasing (see Figure \ref{figV} for $\alpha$ and Section \ref{SAcost}), while the average stock of pollutants in the atmosphere up to the policy adoption is decreasing in Figure \ref{figP}. 
This is related intuitively to the decision maker's willingness to be more proactive and optimally adopt the policy sooner than later, as observed in Figure \ref{figa} (see Section \ref{SAtau} for details), such that emissions are reduced before the pollution stock gets a chance to increase too much.

\section{On the proof of Theorem \ref{thm:mainthm}: Characterising the optimal emissions reduction time}
\label{sec:proof}

In this section we develop a constructive proof of Theorem \ref{thm:mainthm}, leading to the solution to the problem \eqref{valf1} and the characterisation of the optimal emissions reduction time. This will be distilled through a series of subsections and intermediate results.

\subsection{Reformulation of Problem \eqref{valf1}}
\label{section4n}

In the next couple of Sections \ref{never}--\ref{now}, we exploit the Markovian formulation of the problem in \eqref{partprob1} in order to obtain the values of the two most extreme strategies of the social planner, i.e.~never adopt or immediately adopt the environmental policy. 
These, in turn, will be helpful in the dimensionality reduction of our problem in Section \ref{reduce}.

\subsubsection{Never adopt the policy}  
\label{never}

Suppose that the social planner decides to postpone the policy forever, i.e. chooses $\tau=\infty$ in \eqref{partprob1}. 
Then, in view of \eqref{expectation}, using the explicit expression \eqref{po1} of $P^p$, invoking Fubini's theorem, 
the total value $V_\infty$ of this strategy is
\begin{align}
\begin{split}
V_\infty(x,p,\pi) 
&:= \mathbb{E}_{x,p,\pi}\left[\int_0^\infty e^{-rt} X_t P_t dt\right]\\
&=\mathbb{E}_{\pi}\left[ \int_0^\infty e^{-rt}X_t^x\left(\frac{\beta E}{\delta}\left( 1-e^{-\delta t}\right)+pe^{-\delta t}\right)dt\right]\\
&=\int_0^\infty e^{-rt}\,\mathbb{E}_{\pi}\left[X_t^x\right]\left(\frac{\beta E}{\delta}\left( 1-e^{-\delta t}\right)+pe^{-\delta t}\right)dt\\
&=\int_0^\infty e^{-rt} x\left(\pi e^{\alpha t}+(1-\pi)e^{-\alpha t}\right) \left(\frac{\beta E}{\delta}\left( 1-e^{-\delta t}\right)+pe^{-\delta t}\right)dt\\
&=\begin{cases}\label{nevera}
\beta E x (\theta \pi +\rho )
+ xp (\theta_0 \pi + \rho_0) , & \text{if}\,\,r>\alpha, \\
+\infty & \text{if}\,\,r\leq\alpha,
\end{cases}
\end{split}
\end{align}
where $\mathbb{E}_{\pi}$ denotes the expectation conditioned on $\Pi_0=\pi$, and where we define, under the assumption that $r>\alpha$, the constants  
\begin{align} \label{const} 
\begin{split}
\theta &:=\frac{2\alpha(2r+\delta)}{(r-\alpha)(r+\alpha)(r+\delta-\alpha)(r+\delta+\alpha)}>0, 
\qquad 
\theta_0 :=\frac{2\alpha}{(r+\delta-\alpha)(r+\delta+\alpha)}>0, 
\\
\rho &:=\frac{1}{(r+\delta+\alpha)(r+\alpha)}>0, 
\qquad \qquad \qquad \qquad \qquad
\rho_0 :=\frac{1}{r+\delta+\alpha}>0.
\end{split}
\end{align}

\subsubsection{Adopt the policy immediately} \label{now}

Suppose now that the social planner decides to adopt the policy immediately, i.e.~chooses $\tau=0$ in \eqref{partprob1}. 
Then, in view of \eqref{expectation}, using the explicit expression \eqref{po2} of $\hat{P}^p$ and invoking Fubini's theorem, the resulting value $V_0$ is 
\begin{align}
\begin{split}
V_0(x,p,\pi) 
:=
I + \mathbb{E}_{x,p,\pi}\left[\int_0^\infty e^{-rt} X_t \hat{P}_t dt\right] 
&= I + \int_0^\infty e^{-rt} x\left(\pi e^{\alpha t}+(1-\pi)e^{-\alpha t}\right) pe^{-\delta t}dt \\
&=\begin{cases}\label{immeda}
xp (\theta_0 \pi + \rho_0) + I, & \text{if}\,\,r+\delta>\alpha, \\
+\infty & \text{if}\,\,r+\delta\leq\alpha,
\end{cases}
\end{split}
\end{align}
for the previously defined positive constants $\theta_0$ and $\rho_0$ in \eqref{const}.

\subsubsection{Reformulation of Problem \eqref{partprob1}: Dimensionality reduction} 
\label{reduce}

Given that in environmental economics, especially in environmental policy adoption discussions, the aforementioned two strategies in Sections \ref{never}--\ref{now} are clearly plausible, we require that they are also admissible. Hence, we make the following assumption on the problem's parameters, essentially ruling out the possibility of these strategies having value $+\infty$, i.e. yielding an infinite expected cost. 

\begin{assumption}\label{ass}
We assume that the discount rate $r$ and the highest average rate of environmental pollution costs $\alpha$ satisfy $r>\alpha$. 
\end{assumption}

After rewriting the problem in the Markovian framework \eqref{partprob1}, it is obvious that the value function depends on all initial values $p,x>0$ and $\pi\in(0,1)$. Therefore, \eqref{partprob1} seems to be a three-dimensional optimal stopping problem. 
However, thanks to the linearity of the running cost function, we will show that it reduces to a truly two-dimensional one, involving the process $(X^{x,\pi},\Pi^\pi)$, while the deterministic evolution of the pollution stock $P^p$ will eventually affect the optimal timing of adopting the environmental policy only indirectly. 
A similar dimensionality reduction was conjectured in \cite{Pindyck1} for the full information case (from a two- to a one-dimensional problem); here we rigorously prove that such a reduction is possible, while extending it to our setting of the three-dimensional problem \eqref{partprob1}. 

To see this, fix some initial values $(x,p,\pi)\in\mathbb{R}_+\times\mathbb{R}_+ \times (0,1)$, and consider first the expectation 
\begin{align} \label{A}
\mathbb{E}_{x,p,\pi}\left[ \int_0^\tau e^{-rt} X_t P_t dt\right]
&=\mathbb{E}_{x,p,\pi}\left[ \int_0^\infty e^{-rt} X_t P_t dt\right] 
- \mathbb{E}_{x,p,\pi}\left[ \int_\tau^\infty e^{-rt} X_t P_t dt\right] \nonumber\\
&= V_\infty(x,p,\pi) 
- \mathbb{E}_{x,p,\pi}\left[ \int_\tau^\infty e^{-rt} X_t P_t dt\right],
\end{align}
where the latter equality follows from the definition \eqref{nevera} of $V_\infty$.
By invoking the tower property, it follows again from \eqref{nevera} under Assumption \ref{ass} that  
\begin{align}
\mathbb{E}_{x,p,\pi}\left[ \int_\tau^\infty e^{-rt} X_t P_t dt\right] 
&=\mathbb{E}_{x,p,\pi}\left[e^{-r\tau}\left(\,\mathbb{E}_{x,p,\pi}\left[ \int_{\tau}^\infty e^{-r (t-\tau)} X_{t} P_{t} dt \,\bigg|\, \mathcal{F}_\tau\right]\right)\right] \nonumber\\
&=\mathbb{E}_{x,p,\pi}\Big[ e^{-r\tau} V_\infty(X_\tau,P_\tau,\Pi_\tau)\Big], 
\end{align}
where the last equality is due to the strong Markov property of the process $(X,P,\Pi)$. 

Using similar arguments as above, together with \eqref{immeda} under Assumption \ref{ass}, we further obtain that
\begin{equation}\label{B}
\mathbb{E}_{x,p,\pi}\left[ \int_\tau^\infty e^{-rt} X_t \hat{P}_t dt\right] 
= \mathbb{E}_{x,p,\pi}\left[ e^{-r\tau} \big( V_0(X_\tau,\hat{P}_\tau, \Pi_\tau) - I \big)\right].
\end{equation}

Then, combining \eqref{A}--\eqref{B} and recalling from \eqref{p}--\eqref{po2} that $\hat{P}_\tau = {P}_\tau$, we conclude that the value function \eqref{partprob1} can be rewritten as 
\begin{align} \label{V=U}
\begin{split}
V(x,p,\pi) 
&= \inf_{\tau\geq 0}\,\mathbb{E}_{x,p,\pi}\left[ \int_0^\tau e^{-rt} X_t P_t dt + e^{-r\tau}I + \int_\tau^\infty e^{-rt} X_t \hat{P}_t dt\right] \\
&= V_\infty(x,p,\pi) 
+ \inf_{\tau\geq 0}\,\mathbb{E}_{x,\pi}\left[ e^{-r\tau} \big(V_0(X_\tau,{P}_\tau, \Pi_\tau) - V_\infty(X_\tau,P_\tau, \Pi_\tau)\big) \right] \\
&= V_\infty(x,p,\pi) 
+ \inf_{\tau\geq 0}\,\mathbb{E}_{x,\pi}\left[ e^{-r\tau} \big(I - \beta E {X}_\tau (\theta \Pi_\tau +\rho) \big) \right] \\
&= V_\infty(x,p,\pi) 
- \sup_{\tau\geq 0}\,\mathbb{E}_{x,\pi}\left[ e^{-r\tau} \big(\beta E {X}_\tau (\theta \Pi_\tau +\rho) - I \big) \right],
\end{split}
\end{align}
where $\mathbb{E}_{x,\pi} [\,\cdot\,] :=\mathbb{E}[\,\cdot\,\big| \, X_0=x, \Pi_0=\pi]$. 

Hence, the solution to the three-dimensional problem $V$ in \eqref{partprob1} is given in terms of the solution to the two-dimensional optimal stopping problem with value function $U$ given by 
\begin{align}\label{valf2}
U(x,\pi)
&:=\sup_{\tau\geq 0}\,\mathbb{E}_{x,\pi} \left[ e^{-r\tau}G(X_\tau, \Pi_\tau) \right] , 
\end{align}
where the supremum is taken over all ${\mathcal{F}}^X$-stopping times and the function 
\begin{equation} \label{G}
G:\R_+\times(0,1) \mapsto \R \qquad \text{is defined by} \qquad G(x,\pi) := \beta E x (\theta \pi + \rho) - I .
\end{equation}
Therefore, the main aim in the remaining of this paper is to solve \eqref{valf2} and achieve the complete characterisation of the optimal strategy. 

It is well known in optimal stopping theory that multi-dimensional optimal stopping problems cannot be solved in general via the standard \textit{guess-and-verify} approach. 
On the one hand, this solution method involves solving a partial differential equation (PDE) associated with the Hamilton-Jacobi-Bellman equation, but closed form solutions to PDEs can rarely be found. 
On the other hand, it is uncertain whether the usually (a priori) assumed smooth-fit condition of the value function (in two-dimensional stopping problems) holds along the free boundary function which defines the optimal stopping strategy. 
In the sequel, we employ a methodology to solve the problem consisting of a direct probabilistic approach and a transformation of the state-space. 

\begin{remark}
In contrast to the standard full information case \cite{Pindyck1}, it will be shown in the forthcoming analysis that the optimal timing for the policy adoption is not given by a simple constant threshold strategy for the socioeconomic cost process $X$ any more. 
Instead, we will show in Section \ref{sec:Solution} that such a threshold will now be a function of the Bayesian learning process $\Pi$. 
What is even more interesting is that, in the process of completely characterising the optimal emissions reduction policy adoption, we will eventually express the optimal timing solely in terms of the learning process crossing a time- and model parameters-dependent boundary (see Sections \ref{section56}--\ref{sec:Smooth} and Theorem \ref{thm:mainthm}).  
\end{remark}

\subsection{Solution to the two-dimensional optimal stopping problem \eqref{valf2}} 
\label{sec:Solution}

We are now ready to begin the analysis of Problem \eqref{valf2}. To that end, by relying on optimal stopping theory \cite[Section 2.2, Chapter I]{Peskir1}, we firstly introduce the continuation region $\mathcal C$ and stopping region $\mathcal S$ defined by 
\begin{align} \label{CS}
\begin{split}
\mathcal{C}&:=\lbrace (x,\pi)\in \mathbb{R}_+\times (0,1):\,U(x,\pi)>G(x,\pi)\rbrace,\\
\mathcal{S}&:=\lbrace (x,\pi)\in\mathbb{R}_+\times (0,1):\,U(x,\pi)=G(x,\pi)\rbrace,
\end{split}
\end{align}
as well as the stopping time 
\begin{equation} \label{tau*}
\tau^*:=\inf\lbrace t\geq 0:\,({X}^{x,\pi}_t,\Pi^{\pi}_t)\in\mathcal{S}\rbrace,
\end{equation}
with the usual convention $\inf\emptyset=+\infty$.
Later, we will show the optimality of $\tau^*$, as expected.

Before we proceed with addressing the problem, we notice that the boundary points $0$ and $1$ are entrance-not-exit for the diffusion $\Pi^\pi$, that is, $\Pi^\pi$ never atttains $0$ or $1$ in finite time whenever its initial value satisfies $\pi\in (0,1)$ (cf.\ \cite[p.12]{Borodin}). The proof of the following result is omitted as it follows similar arguments to the one of \cite[Lemma 3.1]{Decamps1}. 

\begin{lemma}\label{lemmay}
For all $\pi \in (0,1)$, we have
$$\mathbb{P}(\Pi_t^\pi\in (0,1)\,\, \forall\,\, t\geq 0)=1,  
$$ 
while, for $\pi \in \{0,1\}$, we have 
\begin{equation*}
\mathbb{P}(\Pi^1_t=1,\,\,\forall\,\,t\geq 0)=\mathbb{P}(\Pi^0_t=0,\,\,\forall\,\,t\geq 0)=1.
\end{equation*}
\end{lemma}

Also, we observe from \eqref{XPi} that 
\begin{equation}\label{X}
{X}_t^{x,\pi}= x \exp\left\{\int_0^t \Big(2\alpha\Pi_s^\pi-\alpha-\frac{\sigma^2}{2} \Big) ds+\sigma W_t\right\},\,\,\,\mathbb{P}_{x,\pi} \text{-a.s.}
\end{equation}
In light of Assumption \ref{ass}, we can prove that 
$$
\bigl(e^{-rt}{X}_t^{x,\pi}\bigr)_{t\geq 0} 
\quad \text{is a continuous super-martingale with last element} 
\quad \lim_{t\rightarrow \infty} e^{-rt}{X}_t^{x,\pi} =0 .
$$
This further implies, in view of $\Pi_t^\pi\in (0,1)$ for $\pi\in(0,1)$ due to Lemma \ref{lemmay}, that we can adopt the convention 
\begin{equation} \label{trans}
e^{-r\tau}G({X}^{x,\pi}_\tau,\Pi^{\pi}_\tau)=0 \qquad \text{on} \qquad \lbrace \tau =\infty\rbrace.
\end{equation}

\subsubsection{Well-posedness and initial properties of the value function $U$ defined in \eqref{valf2}}
\label{section54}

The next standard result ensures the well-posedness of the optimal stopping problem under study and its proof can be found in Appendix \ref{app:proofs}.

\begin{lemma}  \label{lem:standass}
The problem \eqref{valf2} is well-posed, the stopping time $\tau^*$ in \eqref{tau*} is optimal and  
\begin{align} \label{Umart}
\begin{cases}
t\mapsto e^{-rt}U({X}^{x,\pi}_t,\Pi_t^\pi) &\text{is a supermartingale,}\\
t\mapsto e^{-r(t\wedge\tau^*)}U({X}^{x,\pi}_{t\wedge\tau^*},\Pi^\pi_{t\wedge\tau^*}) &\text{is a martingale.} 
\end{cases}
\end{align}
\end{lemma}

Next we obtain some further  properties of the value function $U(x,\pi)$ and its proof can also be found in Appendix \ref{app:proofs}. 

\begin{proposition}\label{prop43}
Consider the value function $U$ in \eqref{valf2}. Then, we have:
\begin{enumerate}
\item[\rm (i)] $U(\cdot,\cdot)$ is non-negative on $\R_+ \times (0,1)$;
\item[\rm (ii)] $x\mapsto U(x,\pi)$ is non-decreasing on $\mathbb{R}_+$; 
\item[\rm (iii)] $\pi\mapsto U(x,\pi)$ is non-decreasing on $(0,1)$;
\item[\rm (iv)] $x \mapsto U(x,\pi)$ is Lipschitz continuous on $\R_+$;
\item[\rm (v)] $(x,\pi)\mapsto U(x,\pi)$ is continuous on $\R_+ \times (0,1)$.  
\end{enumerate}
\end{proposition}

\subsubsection{The structure of the state-space and the optimal strategy.}\label{section55}

In this section, we aim at giving a rigorous geometric description of the continuation and stopping regions defined in \eqref{CS}. 
The following lemma is a direct consequence of the continuity of the value function of the optimal stopping problem \eqref{valf2} in Proposition \ref{prop43}.(v).

\begin{lemma}
The continuation (resp., stopping) region $\mathcal{C}$  (resp., $\mathcal{S}$) defined in \eqref{CS} is open (resp., closed). 
\end{lemma}

The next proposition shows that the stopping region $\mathcal{S}$ is up-connected in both arguments $x$ and $\pi$; 
consequently the continuation region $\mathcal{C}$ is down-connected in both $x$ and $\pi$ (see Appendix \ref{app:proofs} for the proof).

\begin{proposition}\label{prop2}
Let $(x_0,\pi_0)\in\mathbb{R}_+\times (0,1)$. The following properties hold:
\begin{enumerate}
\item[\rm (i)] $(x_0,\pi_0)\in\mathcal{S} \quad \Rightarrow \quad (x,\pi_0)\in\mathcal{S} \quad \text{for all}\quad x \in [x_0, \infty)$,
\item[\rm (ii)] $(x_0,\pi_0)\in\mathcal{S} \quad \Rightarrow\quad (x_0,\pi)\in\mathcal{S} \quad \text{for all} \quad \pi \in [\pi_0, 1)$.
\end{enumerate}
\end{proposition} 

In light of Proposition \ref{prop2}, we can define the boundary function $b:(0,1)\rightarrow \mathbb{R}_+$ by 
\begin{equation} \label{b}
b(\pi):=\sup\lbrace x>0:\,U(x,\pi)>G(x,\pi) \rbrace, \qquad \text{for all } \pi\in (0,1), 
\end{equation}
under the convention $\sup \emptyset = 0$. 
Then, using Proposition \ref{prop2}.(i), we can obtain the shape of the continuation and stopping regions from \eqref{CS} in the form  
\begin{align} \label{CSb}
\begin{split}
\mathcal{C}=\lbrace (x,\pi) \in\mathbb{R}_+\times (0,1) \,|\, x < b(\pi)\rbrace ,\\
\mathcal{S}=\lbrace (x,\pi) \in\mathbb{R}_+\times (0,1) \,|\, x\geq b(\pi)\rbrace ,
\end{split}
\end{align}
and the optimal stopping time from \eqref{tau*} takes the form of   
\begin{equation} \label{tau*b}
\tau^*:=\inf\lbrace t\geq 0:\, {X}^{x,\pi}_t \geq b(\Pi^{\pi}_t) \rbrace .
\end{equation}
Given all aforementioned results, we can also prove the following (see Appendix \ref{app:proofs}).

\begin{corollary}\label{cor1}
The boundary function $b$ defined by \eqref{b} satisfies the properties: 
\begin{enumerate}
\item[\rm (i)] $\pi\mapsto b(\pi)$ is non-increasing on $(0,1)$;
\item[\rm (ii)] $\pi\mapsto b(\pi)$ is right-continuous on $(0,1)$.
\end{enumerate}
\end{corollary}

\subsection{A Parabolic Formulation of the two-dimensional optimal stopping problem \eqref{valf2}} \label{section56}

In order to proceed further with our analysis and provide the complete characterisation of the optimal policy adoption timing, it is useful to make a transformation of the state-space.
Notice that the process $({X},\Pi)$ defined in \eqref{XPi} is degenerate, in the sense that both components are driven by the same Brownian motion $W$. 
Therefore, we aim at maintaining the diffusive part in only one of the component processes, while transforming the other component to a completely deterministic bounded variation process, leading to a parabolic formulation. 
Similar transformations were employed in the literature \cite{DeAngelis1, DeAngelis2, FFR, JP17}. 

\subsubsection{The transformed state-space process $(Z,\Pi)$.}

We first define the process
\begin{align} 
{Z}_t^z 
:=\frac{\sigma^2}{2\alpha}\ln \Big(\frac{\Pi_t^\pi}{1-\Pi_t^\pi}\Big) -\ln({X}_t^{x,\pi}),
\qquad Z_0 = z:=\frac{\sigma^2}{2\alpha}\ln\Big(\frac{\pi}{1-\pi}\Big)-\ln(x),
\end{align} 
which evolves deterministically (proof follows via Itô's formula) according to 
\begin{align*}
d{Z}_t=\frac{1}{2}\sigma^2dt , \quad Z_0 = z,
\qquad \text{or, equivalently} \qquad  
{Z}_t^z =z+\frac{1}{2}\sigma^2 t, \quad t \geq 0. 
\end{align*}
Overall, the new state-space process $(Z,\Pi)$ is given by 
\begin{align} \label{ZPi}
\begin{cases}
 dZ_t=\frac 12 \sigma^2 dt, &\quad Z_0 = z :=\frac{\sigma^2}{2\alpha}\ln\big(\frac{\pi}{1-\pi}\big)-\ln(x) \in \R,\\ 
d\Pi_t=\frac{2\alpha}{{\sigma}}\Pi_t (1-\Pi_t) d{W}_t, &\quad \Pi_0 = \pi \in (0,1),
\end{cases}
\end{align} 
and its infinitesimal generator is defined for any $f\in C^{1,2} (\mathbb{R}\times(0,1))$ by 
\begin{equation} \label{genZPi}
(\mathcal{G}f)(z,\pi):=\frac{1}{2}\sigma^2\frac{\partial f}{\partial z}(z,\pi)+\frac{1}{2}\Big(\frac{2\alpha}{\sigma}\Big)^2 \pi^2 (1-\pi)^2 \frac{\partial^2 f}{\partial \pi^2}(z,\pi).
\end{equation}

\subsubsection{The transformed value function $W(z,\pi)$.}

For any $(x,\pi) \in \R_+ \times (0,1)$, define the transformation
\begin{equation} \label{T}
T:=(T_1, T_2):\R_+\times(0,1)\to\R\times(0,1), 
\qquad 
(T_1(x,\pi), T_2(x,\pi)) 
=\Big(\frac{\sigma^2}{2\alpha}\ln\Big(\frac{\pi}{1-\pi}\Big)-\ln(x), \pi \Big),
\end{equation}
which is invertible and its inverse is given by
\begin{equation} \label{T-1}
T^{-1}(z,\pi) 
=\Big( e^{-z}\Big(\frac{\pi}{1-\pi}\Big)^\frac{\sigma^2}{2\alpha}, \pi \Big), \quad (z,\pi)\in \R\times(0,1).
\end{equation}
Using the latter inverse transformation, we firstly introduce the {\it transformed} version $W(z,\pi)$ of the value function $U(x,\pi)$ defined in \eqref{valf2} by 
\begin{align} \label{WU}
W(z,\pi)
&:= U\Big(e^{-z}\Big(\frac{\pi}{1-\pi}\Big)^\frac{\sigma^2}{2\alpha}, \pi \Big) ,
\end{align}
and secondly express the process $X$ as a function of the new state-space process $(Z,\Pi)$; namely,  
\begin{align} \label{transX}
{X}_t^{x,\pi} = e^{-{Z}_t^z} \Big(\frac{\Pi_t^\pi}{1-\Pi_t^\pi}\Big)^\frac{\sigma^2}{2\alpha}, \qquad t \geq 0. 
\end{align}
In view of these, the optimal stopping problem \eqref{valf2} can be rewritten in terms of the process $(Z,\Pi)$ defined in \eqref{ZPi} as   
\begin{align} \label{W}
W(z,\pi)
&= \sup_{\tau\geq 0}\,\mathbb{E}_{z,\pi}\left[ e^{-r\tau} F(Z_\tau, \Pi_\tau) \right], 
\qquad \text{where} \quad 
F(z,\pi):=\beta E e^{-z}\Big(\frac{\pi}{1-\pi}\Big)^\frac{\sigma^2}{2\alpha}(\theta \pi+\rho)-I,
\end{align}
and $\mathbb{E}_{z,\pi}$ denotes the expectation under $\mathbb{P}$, conditional on $Z_0=z$ and $\Pi_0=\pi$.

In view of the relationship \eqref{WU}, the value function $W(\cdot,\cdot)$ inherits important properties which have already been proved for $U(\cdot,\cdot)$ in Section \ref{sec:Solution}. In particular, we have directly from Proposition \ref{prop43}.(v) the following result. 

\begin{proposition}
\label{Wcont}
The transformed value function $(z,\pi) \mapsto W(z,\pi)$ is continuous on $\R\times (0,1)$. 
\end{proposition}

Similarly to Section \ref{sec:Solution}, we can also define the corresponding continuation and stopping regions by
\begin{align} \begin{split} \label{C'S'}
\mathcal{C}'&:=\lbrace (z,\pi)\in\mathbb{R}\times (0,1) \,|\, W(z,\pi)>F(z,\pi)\rbrace,\\
\mathcal{S}'&:=\lbrace (z,\pi)\in\mathbb{R}\times (0,1) \,|\, W(z,\pi)=F(z,\pi)\rbrace. \end{split}
\end{align}
Given Proposition \ref{Wcont}, the continuation region $\mathcal{C}'$ is open and the stopping region $\mathcal{S}'$ is closed and given that $T$ from \eqref{T} is a global diffeomorphism, we actually have  
$$
{\mathcal{C}'}= T({\mathcal{C}}) 
\quad \text{and} \quad 
{\mathcal{S}'}= T({\mathcal{S}}),
$$
where $\mathcal{C}$ and $\mathcal{S}$ are the continuation and stopping regions from \eqref{CS} under $(x,\pi)$-coordinates.
Hence, the corresponding optimal stopping time $\tau^*$ from \eqref{tau*} becomes 
\begin{equation} \label{tau*ZPi}
\tau^*:=\inf\lbrace t\geq 0 \,|\, (Z^{z}_t,\Pi^{\pi}_t)\in\mathcal{S}'\rbrace .
\end{equation}

\subsubsection{The transformed optimal stopping boundary.}
\label{sec:transfbd}

In order to obtain the explicit structure of the regions ${\mathcal{C}'}$ and $\mathcal{S}'$ from \eqref{C'S'}, 
we recall the inverse transformation $T^{-1}$ in \eqref{T-1}, the expression of $\mathcal{C}$ in \eqref{CSb} and the positivity of $b$, to obtain 
\begin{eqnarray*}
(z,\pi)\in{\mathcal{C}'} 
& \Leftrightarrow &\Big( e^{-z}\Big(\frac{\pi}{1-\pi}\Big)^\frac{\sigma^2}{2\alpha}, \pi \Big) \in {\mathcal{C}} 
\quad \Leftrightarrow \quad e^{-z}\Big(\frac{\pi}{1-\pi}\Big)^\frac{\sigma^2}{2\alpha} < b(\pi)  
\quad \Leftrightarrow \quad
z > \ln \bigg( \frac1{b(\pi)}\Big(\frac{\pi}{1-\pi}\Big)^\frac{\sigma^2}{2\alpha} \bigg) .
\end{eqnarray*}
Then, by defining
\begin{equation}
\label{c-1}
c^{-1}(\pi):= \ln \bigg( \frac1{b(\pi)}\Big(\frac{\pi}{1-\pi}\Big)^\frac{\sigma^2}{2\alpha} \bigg)
= \ln \bigg(\Big(\frac{\pi}{1-\pi}\Big)^\frac{\sigma^2}{2\alpha} \bigg) - \ln \big( b(\pi) \big),
\end{equation}
we can obtain the structure of the continuation and stopping regions of $W$, which take the form 
\begin{align} \label{C'S'c-1}
\begin{split} 
\mathcal{C}'&:=\lbrace (z,\pi)\in\mathbb{R}\times (0,1) \,|\, z > c^{-1}(\pi)\rbrace,\\
\mathcal{S}'&:=\lbrace (z,\pi)\in\mathbb{R}\times (0,1) \,|\, z \leq c^{-1}(\pi)\rbrace. \end{split} 
\end{align}

By using the expression \eqref{c-1} of the function $c^{-1}(\cdot)$ and taking into account the fact that $b(\cdot)$ is non-increasing due to Corollary \ref{cor1}.(i), we obtain for any $\varepsilon>0$, that 
\begin{equation*} 
c^{-1}(\pi+\varepsilon) - c^{-1}(\pi)  
= \int_{\pi}^{\pi+\varepsilon} \hspace{-3mm}\frac{\sigma^2}{2\alpha u (1-u)}du 
- \big( \ln b(\pi+\varepsilon) - \ln b(\pi) \big) 
\geq \int_{\pi}^{\pi+\varepsilon} \hspace{-3mm}\frac{\sigma^2}{2\alpha u (1-u)} du > 0.
\end{equation*}
That is, $c^{-1}(\cdot)$ is strictly increasing. 
Moreover, the definition \eqref{c-1} of $c^{-1}(\cdot)$ and the right-continuity of $b(\cdot)$ on $(0,1)$ due to Corollary \ref{cor1}.(i), imply that $c^{-1}(\cdot)$ is right-continuous.
These properties are summarised below.

\begin{lemma} \label{lem:c-1}
The function $c^{-1}(\cdot)$ defined in \eqref{c} is strictly increasing and right-continuous on $(0,1)$.
\end{lemma}

In light of Lemma \ref{lem:c-1}, we may now define the function 
\begin{equation}
\label{c}
c(z):=\inf\{\pi \in (0,1) \,|\, z \leq c^{-1}(\pi)\} .
\end{equation}
In the following result, we prove some properties of $z \mapsto c(z)$ and that it identifies with the optimal stopping boundary of Problem \eqref {W}. 
The proof can be found in Appendix \ref{app:proofs}.

\begin{proposition} \label{prop:c}
The free boundary $c$ defined in \eqref{c} satisfies the following properties: 
\begin{enumerate}[\rm (i)]
\item $c(\cdot)$ is non-decreasing on $\R$; 
\item We have $0 \leq c(z) \leq 1$ for all $z\in\R$ with 
$\lim_{z \downarrow -\infty} c(z)=0$ and 
$\lim_{z \uparrow \infty} c(z)=1$;
\item $c(\cdot)$ is continuous on $\R$.
\item The structure of the continuation and stopping regions for \eqref{W} take the form
\begin{align} \label{C'S'c}
\begin{split} 
\mathcal{C}'&:=\lbrace (z,\pi)\in\mathbb{R}\times (0,1) \,|\, \pi < c(z)\rbrace,\\
\mathcal{S}'&:=\lbrace (z,\pi)\in\mathbb{R}\times (0,1) \,|\, \pi \geq c(z)\rbrace. \end{split}
\end{align}
\end{enumerate}
\end{proposition}

\subsection{Smooth-fit property and integral equation for the transformed stopping boundary}
\label{sec:Smooth}

We firstly define (using Dynkin's formula and standard localisation arguments that make the stochastic integral appearing in the application of It\^o's formula a true martingale; see, e.g., Section 25 in \cite{Peskir1} or the proof of Theorem \ref{intw} in the Appendix) the distance of the transformed value function $W$ from its intrinsic value $F$ by
\begin{align} \label{newval3}
w(z,\pi) &:=W(z,\pi)-F(z,\pi) 
=\sup_{\tau\geq 0}\,\mathbb{E}_{z,\pi} \bigg[ \int_0^\tau e^{-rt} q(Z_t, \Pi_t) dt \bigg],
\end{align}
where the function $q(\cdot,\cdot)$ is defined by
\begin{equation} \label{q}
q(z,\pi) := \beta E e^{-z} \Bigl(\frac{\pi}{1-\pi}\Bigr)^\frac{\sigma^2}{2\alpha} \Big( \big( (\alpha-r)\theta+2\alpha\rho \big) \pi - (\alpha+r)\rho \Big)+rI
\qquad \text{for all } (z,\pi) \in \R  \times (0,1).
\end{equation}
In the following result, we provide properties of $w$ (see Appendix \ref{app:proofs} for their proofs) that will be later used in order to derive its smooth-fit property.

\begin{proposition}\label{prop3}
The function $w$ defined by \eqref{newval3} satisfies the following properties:
\begin{enumerate}
\item[\rm (i)] $z \mapsto w(z,\pi)$ is non-decreasing on $\R$;
\item[\rm (ii)] $\pi \mapsto w(z,\pi)$ is non-increasing on $(0,1)$. 
\item[\rm (iii)]
$w \in C^{1,2} (\mathcal{C}')$ and uniquely solves on any open set $\mathcal{R}$, whose closure is a subset of $\mathcal{C}'$, the PDE
\begin{align} \label{PDEw}
 &(\mathcal{G}-r)m(z,\pi) = -q(z,\pi) , \quad \text{for} \quad (z,\pi)\in\mathcal{R}, \qquad \text{with} \quad
\left. m\right\vert_{\partial\mathcal{R}} =\left. w\right\vert_{\partial\mathcal{R}}.
\end{align}
\end{enumerate}
\end{proposition}

Since the transformed boundary function $c(\cdot)$ is non-decreasing on $\R$ due to Proposition \ref{prop:c}.(i), we observe that the process $(Z^z,\Pi^\pi)$ does not necessarily enter immediately into the stopping region $\mathcal{S}'$ expressed by \eqref{C'S'c}, when started from a point $(z,\pi) \in \partial\mathcal{C}$. 
Hence, a classical proof of the continuity of 
$\pi \mapsto \frac{\partial w}{\partial \pi}(z,\pi)$, for all $z\in\R$, (see \cite{Peskir1} for examples) is not feasible.   
In order to prove the latter result we follow arguments as those in the proof of \cite[Lemma 5.5]{DeAngelis2}; see Appendix \ref{app:proofs} for the detailed technical proof. 

\begin{proposition} \label{smooth}
Consider the function $w$ defined by \eqref{newval3}. 
For each $z\in\R$, we have that $\pi \mapsto \frac{\partial w}{\partial \pi}(z,\pi)$ is continuous on $(0,1)$.  
\end{proposition}

In the sequel, we employ the continuity of $c$ from Proposition \ref{prop:c}.(iii), the regularity and monotonicity of $w$ from Proposition \ref{prop3}, the smooth-fit property from Proposition \ref{smooth} and the fact that the component process $(Z_t)_{t\geq 0}$ is actually a time-variable, to use the local-time-space formula from \cite[Theorem 3.1, Remark 3.2.(2)]{Peskir2005} on $(e^{-r t} w(Z_t, \Pi_t))_{t\geq 0}$ and obtain the following result. 
This technical proof can also be found in Appendix \ref{app:proofs}.

\begin{theorem} \label{intw}
For any $(z,\pi)\in\R\times(0,1)$, the function $w$ defined in \eqref{newval3} can be represented by
\begin{align*} 
w(z,\pi) = \mathbb{E}_{z,\pi} \bigg[ \int_0^{\infty} e^{-rt} q(Z_t, \Pi_t) {\bf 1}_{\{\Pi_t \leq c(Z_t)\}} dt \bigg]
= \int_0^{\infty} e^{-rt} \left( \int_{0}^{c(z+\frac12 \sigma^2 t)} \hspace{-3mm}q\big(z+\tfrac12 \sigma^2 t, \pi'\big) p_t(\pi,\pi') d\pi' \right) dt ,
\end{align*}
where $p_t(\pi,\pi') = \frac{d\mathbb{P}_{\pi}(\Pi^\pi_t \leq \pi')}{d\pi'}$ denotes the transition density of $(\Pi^\pi_t)_{t \geq 0}$.
\end{theorem}

In light of the above integral representation of $w$, we are now finally ready to completely characterise the boundary function $c$, and therefore complete the proof of our main Theorem \ref{thm:mainthm}.
To that end, we define, for each $z\in\R$, 
\begin{equation} \label{m}
m(z) := \inf \{ \pi \in (0,1) \,|\, q(z,\pi) < 0 \},
\end{equation}
which uniquely exists since $q(z,0+) = rI>0$, $q(z,1-) = -\infty$ and $\pi \mapsto q(z,\pi)$ is continuous and decreasing.

By evaluating the integral representation of $w$ in Theorem \ref{intw} at $\pi=c(z)$ for each $z\in\R$ and using the fact that $w(z,c(z))=0$, we obtain that $c$ solves the integral equation (cf.\ \eqref{solc}) 
\begin{equation}
\label{maineq_reminder}
0 = \int_0^{\infty} e^{-rt} \left( \int_{0}^{c(z+\frac12 \sigma^2 t)} \hspace{-3mm}q\big(z+\tfrac12 \sigma^2 t, \pi'\big) p_t(c(z),\pi') d\pi' \right) dt .
\end{equation}
Moreover, by letting 
$$
\mathcal{M} := \{f:\R \mapsto \R \,|\, \text{f is non-decreasing, continuous and satisfies $f(z) \geq m(z)$ for all $z\in\R$}\} ,
$$
the four-step procedure (without additional challenges) developed in \cite{PeskirAm} via the exploitation of the superharmonic property of $W$, can be employed to conclude that the boundary function $c$ defined by \eqref{c} is the unique solution in $\mathcal{M}$. 

Upon collecting all the results developed in this section, the constructive proof of Theorem \ref{thm:mainthm} is therefore complete.

\subsection{Numerical algorithm} 
\label{sectionNumAlgo}

From a numerical point of view, the main challenge consists in solving the functional equation \eqref{maineq_reminder}, which characterises the function $c$ in Theorem \ref{thm:mainthm}. Let us notice that, for $z \in \mathbb{R}$,  
\begin{multline}
\label{maineq_equiv}
\int_0^{\infty} e^{-rt} \left( \int_{0}^{c(z+\frac12 \sigma^2 t)} \hspace{-3mm}q\big(z+\tfrac12 \sigma^2 t, \pi'\big) p_t(c(z),\pi') d\pi' \right) dt \\
= \mathbb{E}_{z,c(z)} \bigg[ \int_0^{\infty} e^{-rt} q(Z_t, \Pi_t) {\bf 1}_{\{\Pi_t \leq c(Z_t)\}} dt \bigg] = \frac{1}{r} \, \mathbb{E}_{z,c(z)} \big[  q(Z_\zeta, \Pi_\zeta) {\bf 1}_{\{\Pi_\zeta \leq c(Z_\zeta)\}} \big] =: \mathcal{E}[c](z),
\end{multline}
where the first equality follows from Theorem \ref{intw} and where $\zeta \sim {\rm Exp}(r)$ is an exponentially distributed random time with mean $r>0$. 
We consider the following numerical scheme:
\begin{equation}
\label{numscheme}
\begin{aligned}
&c^{(0)}(z)=m(z),  \quad z \in \mathbb{R},\\
&c^{(n+1)}(z) = c^{(n)}(z) + \lambda \widetilde{\mathcal{E}}[c^{(n)}](z),\quad z \in \mathbb{R} \,\,\text{and}\,\, n \in \mathbb{N},
\end{aligned}
\end{equation}
where $\widetilde{\mathcal{E}}$ is the Monte Carlo approximation of operator $\mathcal{E}$ in \eqref{maineq_equiv} and $\lambda > 0$ is a constant parameter discussed below. The choice $c^{(0)} \equiv m$ is justified by the fact that $m$ is a lower threshold for $c$ and thus represents a convenient initial step for the iterations. If the functions $c^{(n)}$ converge pointwise to some function $\bar c$, we deduce from \eqref{numscheme} that $\widetilde{\mathcal{E}}[\bar c](z) = 0$ for each $z \in \mathbb{R}$, that is, $\bar c$ solves \eqref{maineq_equiv} and hence $\bar c \equiv c$ (recall from Theorem \ref{thm:mainthm} that the solution to \eqref{maineq_equiv}  exists and is unique). Once function $c$ is computed by the scheme in \eqref{numscheme}, the other functions and variables in the paper can be easily obtained.

Parameter $\lambda$ helps speeding up the convergence. For the numerical tests in this paper, instead of choosing a fixed parameter, we consider $\lambda =\lambda(z)$ (this does not impact the convergence arguments above, as long as $\lambda(z)$ does not depend on $n$). 
In particular, we consider $\lambda(z)=1-m(z)$ when $m(z)>0.5$ and $\lambda(z)=m(z)$ otherwise. By the plot of $m(z)$ in Figure \ref{c>m}, this choice empirically ensures that $\lambda(z)$ is small in the areas where $c^{(n)}$ is expected to approach 0 and 1, thus avoiding excessive oscillations near the boundaries (we recall that $c(z) \in (0,1)$).

For the tests and plots in this paper, we consider the following values (unless specified otherwise) for the model parameters:
\begin{equation}
\sigma = 0.2, \qquad
E = 0.5, \qquad
\beta = 0.4, \qquad
\delta = 0.2, \qquad
I = 10, \qquad
r = 0.1, \qquad
\alpha = 0.05.
\end{equation}

\section{Conclusions}
\label{section6}

This paper focuses on the social planner's option to adopt an environmental policy implying a once-and-for-all reduction in the current emissions. 
In the course of this, we allow for two layers of uncertainty about the future social and economic consequences of the environmental damage, by letting the associated process $X$ fluctuate stochastically and allowing the decision maker to have only partial information about the trend of $X$. 
Introducing partial information about key parameters in the stochastic dynamics of socioeconomic costs of pollution reflects better the current real-life debate on the actions required to contrast climate change. 
The rigorous and complete treatment of this nontrivial novel problem represents the main contribution of our work. 
We suitably tackle this increased uncertainty on the after-effects of pollution and show that it plays a considerably important role in the optimal timing of policy adoption. A first important effect of additional uncertainty is indeed that the optimal policy adoption is no longer triggered by constant thresholds, as in the full information case \cite{Pindyck1}.

As a matter of fact, without relying on the traditional \textit{guess-and-verify-approach} (non-feasible in our multi-dimensional case), we show via probabilistic means and state-space transformation 
techniques that it is optimal to abate emissions when the stochastic socioeconomic costs $X$ hits or exceeds an upper tolerance level $b(\Pi)$, depending on the filtering estimate process $\Pi$ of the unknown drift of $X$. Such a stopping rule can be equivalently expressed in terms of a boundary $c$, depending on a scaled time coordinate $Z$: It is optimal to reduce the stock of pollutants when the belief process $\Pi$ -- likelihood of having on average increasing socioeconomic costs -- becomes ``decisive enough'' and exceeds the time-dependent percentage $c(Z)$. Sensitivity of the optimal policy adoption with respect to the model's parameters is studied by solving numerically the nonlinear integral equation which is uniquely solved by the boundary $c$. 

Interestingly, our numerical study reveals that an increase in the amplitude of the two layers of uncertainty induces different effects on the expected optimal timing of policy adoption. On the one hand, an increase in the volatility $\sigma$ of the socioeconomic cost process $X$ gives rise to a sharper estimate of the true trend of $X$ (as the volatility of the learning process decreases), so that the first layer of uncertainty (the Brownian risk) prevails and the expected optimal time of pollution reduction increases. This effect is in line with the classical ``value of waiting'' paradigm in real options. On the other hand, by increasing the average rate $\alpha$ of increase/decrease of future socioeconomic costs $X$, the variances of both the learning process $\Pi$ and of the unknown trend of $X$ increase, with the effect of making the decision maker more proactive willing to bring forward the optimal time of pollution reduction (on average).

\subsection{Problems with a similar structure}
\label{Sec:similar}

The solution method presented in this paper applies also to different variations of the environmental policy adoption problem, as well as other types of problems in decision making under two layers of uncertainty. 

\vspace{1mm}
{\it Example 1.}
Firstly, in the current setup, when the social planner has only partial information on the dynamics \eqref{XX} of socioeconomic costs  $X=(X_t)_{t\geq 0}$ generated by a unit of pollution, we can consider a stochastically evolving stock of pollutants (instead of \eqref{P})
\begin{align} \label{stochP}
\begin{split}
dP_t &=(\beta E-\delta P_t)dt + \eta d{\widetilde B}_t, \qquad \text{for all } t \leq \tau, \qquad P_0=p>0, \\
d\hat{P}_t &= -\delta \hat{P}_t dt + \eta d{\widetilde B}_t, \qquad\qquad\;\; \text{for all } t > \tau, \qquad  \hat{P}_\tau=P_\tau .
\end{split}
\end{align}
where the parameters $E, \beta, \delta$ are as in \eqref{P}, while $(\widetilde B_t)_{t\geq 0}$ is a Brownian motion (independent of $B$) modelling the shocks affecting the atmospheric stock of pollutants and the volatility $\eta>0$ denotes their extend.
Such dynamics will neither interfere with the learning process of the decision maker nor affect the analysis resulting to the two-dimensional problem \eqref{valf2}. 
The only difference is that the expected values of the stock of pollutants 
\begin{align} \label{EstochP}
\begin{split}
\E^{\P_\pi}[P_t] &= e^{-\delta t} p +  \frac{\beta_1 E \pi + \beta_2 E (1-\pi)}{\delta} (1- e^{-\delta t}), \qquad \text{for all } t \leq \tau , \\
\E^{\P_\pi}[\hat{P}_t] &= e^{-\delta t} p , \qquad\qquad\qquad\qquad\qquad\qquad\qquad\qquad\; \text{for all } t > \tau ,
\end{split}
\end{align}
should be used in the calculations, instead of the explicit expressions \eqref{po1} and \eqref{po2}. 
All subsequent analysis of the resulting two-dimensional problem should be identical. 

\vspace{1mm}
{\it Example 2.}
Another alternative would be to consider the setup when $X_t$ is the price per unit of a given good at time $t \geq 0$, that evolves according to the dynamics 
$$
dX_t^{x}= \mu X_t^x dt+\sigma X_t^x d{B}_t, \quad X_0^{x} = x > 0, 
$$
where the future price trend $\mu$ is random and unknown to the decision maker. 
The latter can control the production process $P$ of this good, that evolves according to the dynamics
\begin{equation*} 
dP_t= \beta (\lambda_0 + \lambda_1 \textbf{1}_{\{t \geq \tau\}} - P_t)dt + \eta d{\widetilde B}_t, \qquad \text{for all } t \geq 0, \qquad P_0=p>0,
\end{equation*}
such that the mean production rate $\lambda_0$ can be increased to $\lambda_0+\lambda_1$ at time $\tau$, chosen by the decision maker, at a cost $I(X_\tau)$ that could depend on the current price $X_\tau$ of the good in the economy at that time. 
The decision maker's aim would then be to find the optimal timing for production expansion, while learning the trend of this good's price via a learning process $\Pi$, in order to maximise their profits (net of production expansion costs), which takes the form of 
$$
\sup_{\tau\geq 0}\,\mathbb{E}_{x,p,\pi}\Big[\int_0^\tau e^{-rt} {X}_t P_t dt - e^{-r\tau} I(X_\tau) + \int_\tau^\infty e^{-rt} {X}_t \hat{P}_t dt\Big]
$$
This problem can be studied using a similar analysis as in this paper. 

\vspace{1mm}
{\it Example 3.} 
More generally, consider a two-dimensional stochastic processes $(X,P)$, whose components are independent and $X$ is always positive (e.g. generalised geometric Brownian motion), while $P$ is a mean-reverting process. 
The decision maker is faced with a second layer of uncertainty in the drift of $X$ and uses a learning process $\Pi$ to learn the non-observable random drift via observations.  
The decision maker's control is a stopping time that either  changes (upwards/downwards) the drift of $P$ once and for all, if it is part of the running cost/reward in the problem (i.e. $i=1$ below), with a sunk cost/reward that can depend on $X$, or simply stops the process $X$, if it is the only process involved (i.e. $i=0$ below).
The decision maker's aim is to maximise/minimise an optimisation criterion, given any initial values $(x,p,\pi)\in\mathbb{R}_+\times \mathbb{R} \times (0,1)$ and constants $\gamma \in \mathbb{R}$ and $i\in\{0,1\}$, in the form of the optimal stopping problems
\begin{align} \label{otherV}
\begin{split}
V^{\pm}_1(x,p,\pi) &:=\inf_{\tau\geq 0}\,\mathbb{E}_{x,p,\pi}\Big[\int_0^\tau e^{-rt} \gamma {X}_t P_t^{i} dt \pm e^{-r\tau} I(X_\tau) + \int_\tau^\infty e^{-rt} \gamma {X}_t \hat{P}^{i}_t dt\Big], \\
V^{\pm}_2(x,p,\pi) &:=\sup_{\tau\geq 0}\,\mathbb{E}_{x,p,\pi}\Big[\int_0^\tau e^{-rt} \gamma {X}_t P_t^{i} dt \pm e^{-r\tau} I(X_\tau) + \int_\tau^\infty e^{-rt} \gamma {X}_t \hat{P}_t^{i} dt\Big], 
\end{split}
\end{align}
which can be studied using similar techniques as in this paper.

\subsection{Ideas for extensions}

Our work provides a first tractable example of a stylised model aiming at describing the role of different sources of uncertainty in the typically irreversible social decisions related to climate policy. Clearly, there are various direction towards our analysis can be generalised. 

First of all, an equally significant uncertainty could be considered in the drift of the dynamics \eqref{stochP} of the stock of pollution $P$, instead of its socioeconomic costs $X$. Such a problem can also be reduced to a two-dimensional one, following similar arguments as in Section \ref{section4n}, but the resulting problem will involve a significantly more complicated structure for the reward function $G$ (cf.\ \eqref{G}) and the underlying dynamics of $X$ and $\Pi$ will be driven by two distinct and correlated Brownian motions (instead of a common one, cf.\ \eqref{XPi}).

Second of all, it would be interesting to account for Knightian uncertainty in the dynamics of the cost process (see, e.g., \cite{Ferrarietal}, \cite{NishimuraOzaki}, \cite{HellmannThijssen} for examples of real-options/irreversible investment problems under Knightian uncertainty). This would in turn allow to comparatively study how the different specifications of uncertainty (no uncertainty, Bayesian uncertainty, Knightian uncertainty) affect the optimal emissions reduction policy. 

Finally, it would be important to allow for the transboundary effect of pollution and thus incorporate strategic interaction between different decision makers (see the recent \cite{Boucekkine-etal} for a spatial deterministic game of pollution control).
How the free-riding effect will depend on the specification of uncertainty would represent a key question in that context (see \cite{Kwon1}, \cite{Kwon2} for recent works on public good contribution games under uncertainty).

\appendix

\section{Technical Proofs}
\label{app:proofs}

{\bf Proof of Lemma \ref{lem:standass}.}
We show below that the process $({X}_t,\Pi_t)$ from \eqref{XPi} satisfies the condition  
\begin{equation}\label{standass}
\mathbb{E}_{x,\pi}\left[ \sup_{t\geq 0} e^{-rt}G({X}_t,\Pi_t)\right]<+\infty,
\end{equation}
which will then imply, according to \cite[Theorem D.12]{Karatzas1}, that the stopping time \eqref{tau*} is optimal for the well-posed problem \eqref{valf2}, as well as ensures that \eqref{Umart} holds true.

To that end, notice that using the expression of $G$ in \eqref{G}, we obtain
\begin{align*}
\mathbb{E}_{x,\pi}\left[ \sup_{t\geq 0} e^{-r t} \big(\beta E {X}_t (\theta \Pi_t +\rho) - I \big) \right]
&\leq \mathbb{E}_{x,\pi} \left[ \sup_{t\geq 0} e^{-rt} \big(\beta E {X}_t (\theta \Pi_t +\rho)
\big)\right]\\
&\leq \mathbb{E}_{x,\pi} \left[ \sup_{t\geq 0} e^{-rt} \beta E (\theta+\rho) {X}_t\right] ,
\end{align*}
where the latter inequality is due to the positivity of $e^{-rt}\beta E \theta {X}_t$ and the property $\Pi_t^\pi\in (0,1)$ for $\pi\in(0,1)$ (Cf.\ Lemma \ref{lemmay}). 
Thus, using the explicit expression of $X$ in \eqref{X}, the upper bound becomes  
\begin{align*}
\mathbb{E} \left[ \sup_{t\geq 0} e^{-rt} \beta E (\theta+\rho){X}_t^{x,\pi} \right]
&=\beta E (\theta+\rho) x \, 
\mathbb{E}\left[ \sup_{t\geq 0} e^{-rt} \, e^{\int_0^t (2\alpha\Pi_s^\pi-\alpha-\frac{\sigma^2}{2}) ds+\sigma W_t}\right]\\
&\leq \beta E(\theta+\rho)x\,\mathbb{E}\left[ \sup_{t\geq 0} e^{- b t + \sigma W_t}\right],
\end{align*}
where we define the constant $b:=r-\alpha+\frac{\sigma^2}{2} > 0$. 
Hence, using \cite[Section 3.5.C]{Karatzas2}, which implies that 
\begin{equation*}
\text{for $z>0$,} \qquad
\mathbb{P}\Bigl(\sup_{t\geq 0} \big\{-bt+ {\sigma}W_t \big\}\in dz \Bigr) 
=\tfrac{2 b}{{\sigma}^2} e^{ -\frac{2 b}{{\sigma}^2}z } dz,
\end{equation*}
we can conclude that
\begin{align*}
&\mathbb{E}\left[ \sup_{t\geq 0}e^{- b t+{\sigma}W_t}\right]
=\frac{2 b}{{\sigma}^2} \int_0^\infty e^z e^{-\frac{2 b}{{\sigma}^2}z}dz 
=\frac{2 b}{{\sigma}^2} \int_0^\infty e^{-\big(\frac{2 b}{{\sigma}^2}-1\big)}dz<+\infty,
\end{align*}
where we used the fact that ${2 b}>{{\sigma}^2}$ due to Assumption \ref{ass}. 
Hence, we conclude that \eqref{standass} holds.
\hfill \cvd

{\bf Proof of Proposition \ref{prop43}.}
{\it Proof of part }(i).
This is a trivial consequence of the definition \eqref{valf2} and the property \eqref{trans}, which imply the non-negativity of $U(x,\pi)$, since never stopping, i.e.~choosing $\tau=\infty$, is an admissible strategy which results in a payoff of zero. 

{\it Proof of part }(ii).
By the definition \eqref{valf2} of $U$, 
the explicit expression \eqref{X} of ${X}^{x,\pi}$ which implies that $x\mapsto {X}^{x,\pi}_\tau$ is increasing for any stopping time $\tau$, 
and the definition \eqref{G} of $G$ which implies that $x \mapsto G(x,\pi)$ is increasing for any $\pi \in(0,1)$, 
we conclude that $x\mapsto U(x,\pi)$ is increasing as well. 

{\it Proof of part }(iii).
Using the above properties together with the fact that $\pi\mapsto \Pi^{\pi}_\tau$ is increasing for any stopping time $\tau$ (see the comparison theorem of Yamada and Watanabe, e.g. \cite[Proposition 2.18]{Karatzas2}), 
the explicit expression \eqref{X} of ${X}^{x,\pi}$ which then implies that $\pi\mapsto {X}^{x,\pi}_\tau$ is also increasing for any stopping time $\tau$, 
and the definition \eqref{G} of $G$ which implies that $\pi \mapsto G(x,\pi)$ is increasing for any $\pi \in(0,1)$, 
we conclude by the definition \eqref{valf2} of $U$ that $\pi \mapsto U(x,\pi)$ is increasing as well.
 
{\it Proof of part }(iv).
For any $0<x_1<x_2$ and $\pi\in(0,1)$, we define by $\tau^*:=\tau^*(x_2,\pi)$ the optimal stopping time for the value function $U(x_2,\pi)$ in \eqref{valf2}. 
Then using the monotonicity of $U(\cdot,\pi)$ from part (ii) and the expression \eqref{valf2} of $U$, we obtain for a sufficiently large constant $C_1(x,\pi_2)$ (cf.\ Lemma \ref{lem:standass})
\begin{align}
0
&\leq U(x_2,\pi)-U(x_1,\pi) \nonumber\\
&\leq \mathbb{E}_{x_2,\pi}\left[ e^{-r\tau^*}
\big(\beta E {X}_{\tau^*} (\theta \Pi_{\tau^*} +\rho) - I \big) \right] 
- \mathbb{E}_{x_1,\pi}\left[ e^{-r\tau^*}
\big(\beta E {X}_{\tau^*} (\theta \Pi_{\tau^*} +\rho) - I \big) \right] \nonumber\\
&=(x_2-x_1)\,\mathbb{E} \left[ e^{-r\tau^*}
\beta E {X}^{1,\pi}_{\tau^*} (\theta \Pi^\pi_{\tau^*} +\rho) \right] \nonumber\\
&\leq (x_2-x_1) \,\mathbb{E} \left[ e^{-r\tau^*} \beta E {X}^{1,1}_{\tau^*} (\theta +\rho) \right] 
\leq C_1(x,\pi_2)(x_2-x_1), \label{estUQx}
\end{align}
where the penultimate inequality follows from the positivity of coefficients, $\Pi_t \in (0,1)$ for all $t\geq 0$ due to Lemma \ref{lemmay}, and ${X}^{1,\pi}_{\tau^*} \leq {X}^{1,1}_{\tau^*}$, since $\pi\mapsto {X}^{x,\pi}_{\tau^*}$ is increasing for the fixed $\tau^*$. 

{\it Proof of part }(v). 
To show this, it is sufficient to use part (iv) and additionally prove that
\begin{align*}
\pi \mapsto U(x,\pi) \quad \text{is continuous.} 
\end{align*}
To that end, fix $x>0$, $0<\pi_1<\pi_2<1$ and define by $\bar\tau^*:=\tau^*(x,\pi_2)$ the optimal stopping time for $U(x,\pi_2)$ in \eqref{valf2}. 
Using the monotonicity of $U(x,\cdot)$ from part (iii) and the expression in \eqref{valf2}, we deduce that
\begin{align}
0&\leq U(x,\pi_2) - U(x,\pi_1) \nonumber\\
&\leq \mathbb{E}_{x, \pi_2}\left[ e^{-r\bar\tau^*}
\big(\beta E {X}_{\bar\tau^*} (\theta \Pi_{\bar\tau^*} +\rho) - I \big) \right] 
- \mathbb{E}_{x, \pi_1}\left[ e^{-r\bar\tau^*}
\big(\beta E {X}_{\bar\tau^*} (\theta \Pi_{\bar\tau^*} +\rho) - I \big) \right] \nonumber\\
&=\mathbb{E} \left[ e^{-r\bar\tau^*} \beta E  
\big({X}^{x,\pi_2}_{\bar\tau^*} (\theta \Pi^{\pi_2}_{\bar\tau^*} +\rho) - {X}^{x,\pi_1}_{\bar\tau^*} (\theta \Pi^{\pi_1}_{\bar\tau^*} +\rho) \big) \right]. 
\label{estimatepi}
\end{align}

We now aim at taking limits as $\pi_1 \to \pi_2$. To that end, we notice that, due to the positivity of coefficients, $\Pi_t \in (0,1)$ for all $t\geq 0$ thanks to Lemma \ref{lemmay}, and ${X}^{x,0}_{\bar\tau^*} \leq {X}^{x,\pi_1}_{\bar\tau^*} \leq {X}^{x,\pi_2}_{\bar\tau^*} \leq {X}^{x,1}_{\bar\tau^*}$ (since $\pi \mapsto {X}^{x,\pi}_{\bar\tau^*}$ is increasing for the fixed $\bar\tau^*$), we conclude that 
\begin{align*} 
e^{-r\bar\tau^*} \beta E \big({X}^{x,\pi_2}_{\bar\tau^*} (\theta \Pi^{\pi_2}_{\bar\tau^*} +\rho) - {X}^{x,\pi_1}_{\bar\tau^*} (\theta \Pi^{\pi_1}_{\bar\tau^*} +\rho) \big)  
\leq e^{-r\bar\tau^*} \beta E \big({X}^{x,1}_{\bar\tau^*} (\theta +\rho) - {X}^{x,0}_{\bar\tau^*} \rho \big) , 
\end{align*}
such that
\begin{align*} 
\mathbb{E} \left[ e^{-r\bar\tau^*} \beta E \big({X}^{x,\pi_2}_{\bar\tau^*} (\theta \Pi^{\pi_2}_{\bar\tau^*} +\rho) - {X}^{x,\pi_1}_{\bar\tau^*} (\theta \Pi^{\pi_1}_{\bar\tau^*} +\rho) \big) \right] 
\leq x \, C_2(x,\pi_2).
\end{align*}
Here, $C_2(x,\pi_2)>0$ is a sufficiently large constant (note that the calculations in the proof of Lemma \ref{lem:standass} imply the finiteness of the expectation).
It thus follows that the dominated convergence theorem can be applied when letting $\pi_1 \to \pi_2$ in \eqref{estimatepi}, so that that the upper bound in \eqref{estimatepi} tends to zero. This implies the continuity of $\pi \mapsto U(x,\pi)$ and completes the proof.
\hfill \cvd

{\bf Proof of Proposition \ref{prop2}.} 
In order to prove these results, we firstly define the distance $u$ of the value function from its intrinsic value, whose expression is obtained by an application of Dynkin's formula:  
\begin{align*}
\begin{split}
u(x,\pi)&:=U(x,\pi)- G(x,\pi) 
=\sup_{\tau\geq 0}\,\mathbb{E}_{x,\pi}\left[ \int_0^\tau e^{-rt}\Big( \beta E\bigl(\theta(\alpha-r)+2\alpha\rho\bigr) {X}_t\Pi_t - (\alpha+r)\beta E\rho{X}_t + rI \Big)dt\right].
\end{split} 
\end{align*}

{\it Proof of part }(i).
Thanks to Assumption \ref{ass}, we can easily verify that $\theta(\alpha-r)+2\alpha\rho<0$. 
Then, from the explicit solution of ${X}^{x,\pi}$ given in \eqref{X}, it is clear that $x\mapsto {X}^{x,\pi}$ is increasing. This implies that for any stopping time $\tau$, we have for any $\pi \in (0,1)$, that
\begin{equation*}
x\mapsto \mathbb{E}_{x,\pi}\left[ \int_0^\tau e^{-rt}\left( \beta E \bigl(\theta(\alpha-r)+2\alpha\rho\bigr) {X}_t\Pi_t -(\alpha+r)\beta E\rho {X}_t + rI \right)dt\right] 
\quad \text{is decreasing on } \R,
\end{equation*}
which implies that $x \mapsto u(x,\pi)$ is also decreasing.

Suppose now that $(x_0,\pi_0)\in\mathcal{S}$ and consider $(x,\pi_0)$ for some $x\geq x_0$. By the monotonicity of $u(\cdot,\pi)$, we have 
$u(x, \pi_0)\leq u(x_0, \pi_0)=0.$
Since $u(x,\pi_0)$ is non-negative by definition, we must have $u(x,\pi_0)=0$, thus $U(x,\pi)=G(x,\pi)$, i.e. $(x,\pi_0)\in\mathcal{S}$. 

{\it Proof of part }(ii). 
Recall that $\pi \mapsto \Pi^\pi$ is increasing by the comparison theorem of Yamada and Watanabe (see, e.g., \cite[Proposition 2.18]{Karatzas2}) and consequently that $\pi \mapsto {X}^{x,\pi}$ is also increasing in light of \eqref{X}. 
Therefore, for any stopping time $\tau$, we clearly have for any $x>0$, that
\begin{equation*}
\pi \mapsto \mathbb{E}_{x,\pi}\left[ \int_0^\tau e^{-rt}\left( \beta E \bigl(\theta(\alpha-r)+2\alpha\rho\bigr) {X}_t\Pi_t  -(\alpha+r)\beta E\rho {X}_t+rI\right)dt\right] 
\quad \text{is decreasing on } (0,1),
\end{equation*}
implying that $\pi \mapsto u(x,\pi)$ is also decreasing. 

This in turn implies that for $(x_0,\pi_0) \in \mathcal{S}$ and $(x_0,\pi)$ for some $\pi\geq \pi_0$,  we have $u(x_0,\pi)\leq u(x_0,\pi_0)=0$. 
We then conclude that $u(x_0,\pi)=0$, i.e. $(x_0,\pi)\in\mathcal{S}$, which completes the proof. 
\hfill \cvd

{\bf Proof of Corollary \ref{cor1}.}
{\it Proof of part }(i).
This follows directly from the definition \eqref{b} of $b$, the shape of continuation and stopping regions in \eqref{CSb} and Proposition \ref{prop2}.(ii). 

{\it Proof of part }(ii).
Let $(\pi_n)_{n\in\mathbb{N}}$ be a decreasing sequence in $(0,1)$ that converges to some $\pi_0\in (0,1)$. 
Since $b(\cdot)$ is non-increasing, we have that $b(\pi_n)$ is non-decreasing in $n\in \mathbb{N}$ and bounded above by $b(\pi_0)$. 
Thus, the limit  $b(\pi_0+):=\lim_{n\rightarrow \infty} b(\pi_n)$ exists. 

Since $(b(\pi_n),\pi_n) \in \mathcal{S}$ we have $U(b(\pi_n),\pi_n)=G(b(\pi_n),\pi_n)$ and by the continuity of the value function $U$ in Proposition \ref{prop43}.(iv) and $G$ by definition \eqref{G}, we conclude that 
\begin{equation*}
U(b(\pi_0+),\pi_0)=G(b(\pi_0+),\pi_0).
\end{equation*}
This implies that $b(\pi_0+)\geq b(\pi_0)$ and due to the fact that $b(\cdot)$ is non-increasing, we obtain $b(\pi_0+) = b(\pi_0)$ which completes the proof. 
\hfill \cvd

{\bf Proof of Proposition \ref{prop:c}.}
{\it Proof of part} (i). 
This claim follows from Lemma \ref{lem:c-1}, together with the definition \eqref{c} of $c$. 

{\it Proof of part} (ii).
We observe from the definition \eqref{c-1} of $c^{-1}(\cdot)$ that (since $b$ is bounded by the constant thresholds associated to the the full information case with trend $\pm \alpha$) we have 
$$\lim_{\pi \downarrow 0} c^{-1}(\pi) = -\infty 
\qquad \text{and} \qquad 
\lim_{\pi \uparrow 1} c^{-1}(\pi) = \infty ,$$ 
Taking these into account together with the definition \eqref{c} of $c(\cdot)$, we then conclude that 
$$\lim_{z \downarrow -\infty} c(z) = 0 
\qquad \text{and} \qquad 
\lim_{z \uparrow \infty} c(z) = 1 .$$
The non-decreasing property of $c(\cdot)$ from part (i) then completes the proof of this part. 

{\it Proof of part} (iii). This follows from \cite[Proposition 1.(7)]{Embrechts}, upon using the strictly increasing property of $c^{-1}(\cdot)$ (cf.\ Lemma \ref{lem:c-1}) and \eqref{c}.

{\it Proof of part} (iv). 
This claim again follows from the definition \eqref{c} of $c$ and its monotonicity from part (i), combined with the expressions of the sets in \eqref{C'S'c-1}.
\hfill \cvd

{\bf Proof of Proposition \ref{prop3}.}
In both parts, we use the fact that $\theta(\alpha-r)+2\alpha\rho<0$ thanks to Assumption \ref{ass}. 

{\it Proof of part }(i).
This follows immediately due to the fact that 
$$z\mapsto \exp\{-{Z}_t^z\} = \exp\Big\{-z-\frac{1}{2}\sigma^2 t \Big\} \quad \text{is non-negative and decreasing on } \R,$$
which implies that $z \mapsto w(z,\pi)$ is non-decreasing on $\R$.

{\it Proof of part }(ii).
We firstly recall that $\pi \mapsto \Pi^\pi$ is increasing and observe that this yields  
$$\pi\mapsto \Big(\frac{\Pi_t^\pi}{1-\Pi_t^\pi} \Big)^\frac{\sigma^2}{2\alpha} \quad \text{is increasing on } (0,1), \quad \mathbb{P}-\text{a.s.},$$ which implies that 
$\pi \mapsto w(z,\pi)$ is non-increasing on $(0,1)$.

{\it Proof of part }(iii).
In view of \eqref{Umart} together with \eqref{WU}, we know that 
\begin{equation*}
t\mapsto e^{-r(t\wedge\tau^*)} W({Z}_{t\wedge\tau^*},\Pi_{t\wedge\tau^*}) \quad 
\text{is a martingale.}
\end{equation*}
Taking this into account together with the problem's parabolic formulation (cf.\ \eqref{genZPi}), we can make use of standard arguments in the general theory of optimal stopping (see, e.g. \cite[Section 7.1]{Peskir1}, among others) and classical PDE results on the regularity for solutions of parabolic differential equations (cf. \cite[Corollary 2.4.3]{Krylov}) 
to conclude that $W$ is the unique classical solution, on any open set $\mathcal{R}$ whose closure is contained in $\mathcal{C}'$, of the PDE 
\begin{align*} 
\begin{split}
&(\mathcal{G}-r)m(z,\pi) =0 , \quad \text{for} \quad (z,\pi)\in\mathcal{R}, 
\qquad \text{with} \quad
\left. m\right\vert_{\partial\mathcal{R}} =\left. W\right\vert_{\partial\mathcal{R}}.
\end{split}
\end{align*}
In view of this result, the arbitrariness of $\mathcal{R}$, the definition \eqref{newval3} of $w$, and the smooth expression of $F$ in \eqref{W}, we can conclude that ${w}\in C^{1,2}(\mathcal{C}')$ and solves the claimed PDE.
\hfill \cvd

{\bf Proof of Proposition \ref{smooth}.}
The fact that $\pi \mapsto \frac{\partial w}{\partial \pi}(z,\pi)$ is continuous separately in the continuation region $\mathcal{C}'$ and in the stopping region $\mathcal{S}'$ is due to Proposition \ref{prop3}.(iii) and to $w\equiv 0$ on $\mathcal{S}'$, respectively. 
Hence, it remains only to prove the continuity of $\pi \mapsto \frac{\partial w}{\partial \pi}(z,\pi)$ for all $(z,\pi) \in \partial\mathcal{C}'$. 
This is accomplished in the following steps.

{\it Step 1: For any $(z,\pi)\in (\R\times(0,1))\setminus \partial\mathcal{C}'$ and $\tau^* = \tau^*(z,\pi)$ given by \eqref{tau*ZPi}, we have 
$\frac{\partial W}{\partial z}(z,\pi) = - W(z,\pi) - \mathbb{E}_{z,\pi}\left[ e^{-r\tau^*} I \right]$.}
To prove this, we obtain the expression of $\frac{\partial W}{\partial z}$ in two separate parts of the state space: 

For $(z,\pi)\in\mathcal{S}'$, the definition \eqref{C'S'} of $\mathcal{S}'$ implies that $W(z,\pi) = F(z,\pi)$, \eqref{tau*ZPi} implies that $\tau^*=0$ and in view of the definition \eqref{W} of $F$, we have  
$$
\frac{\partial W}{\partial z}(z,\pi) = \frac{\partial F}{\partial z}(z,\pi) 
= - F(z,\pi) - I
= - W(z,\pi) - \mathbb{E}_{z,\pi}\left[ e^{-r\tau^*} I \right].
$$

For $(z,\pi)\in\mathcal{C}'$, we choose a sufficiently small $\varepsilon>0$ such that $(z-\varepsilon,\pi)\in\mathcal{C}'$ and $(z+\varepsilon,\pi)\in\mathcal{C}'$.
Then, we have from \eqref{W} that
\begin{align}
\frac{W(z+\varepsilon,\pi) - W(z,\pi)}\varepsilon 
&\geq \frac1\varepsilon \mathbb{E}\left[ e^{-r\tau^*}
\big( F(Z_{\tau^*}^{z+\varepsilon}, \Pi_{\tau^*}^{\pi}) - F(Z_{\tau^*}^{z}, \Pi_{\tau^*}^{\pi})
\big) \right] \nonumber \\
&= \mathbb{E}\bigg[ e^{-r\tau^*} \beta E e^{- \frac12 \sigma^2 \tau^*} \Big(\frac{\Pi_{\tau^*}^{\pi}}{1-\Pi_{\tau^*}^{\pi}}\Big)^\frac{\sigma^2}{2\alpha}(\theta \Pi_{\tau^*}^{\pi}+\rho) \Big( \frac{e^{-(z+\varepsilon)} - e^{-z}}{\varepsilon} \Big) 
\bigg] \nonumber 
\end{align}  
Similarly, we also obtain 
\begin{align}
\frac{W(z,\pi) - W(z-\varepsilon,\pi)}\varepsilon 
&\leq \mathbb{E}\bigg[ e^{-r\tau^*} \beta E e^{- \frac12 \sigma^2 \tau^*} \Big(\frac{\Pi_{\tau^*}^{\pi}}{1-\Pi_{\tau^*}^{\pi}}\Big)^\frac{\sigma^2}{2\alpha}(\theta \Pi_{\tau^*}^{\pi}+\rho) \Big( \frac{e^{-z} - e^{-(z-\varepsilon)}}{\varepsilon} \Big) 
\bigg]. \nonumber 
\end{align} 
Letting $\varepsilon \downarrow 0$ in both expressions and recalling that $W \in C^{1,2} (\mathcal{C}')$, we 
find that 
\begin{align}
\frac{\partial W}{\partial z}(z,\pi) = - \mathbb{E}\bigg[ e^{-r\tau^*} \beta E e^{- Z_{\tau^*}^{z}} \Big(\frac{\Pi_{\tau^*}^{\pi}}{1-\Pi_{\tau^*}^{\pi}}\Big)^\frac{\sigma^2}{2\alpha}(\theta \Pi_{\tau^*}^{\pi}+\rho)\bigg] 
= - W(z,\pi) - \mathbb{E}_{z,\pi}\left[ e^{-r\tau^*} I \right] .\nonumber 
\end{align} 
\vspace{0.15cm}
 
{\it Step 2: $\frac{\partial W}{\partial z}$ is locally bounded.}
Notice from the expression of $\frac{\partial W}{\partial z}$ in {\it step 1} and the continuity of $W(\cdot,\cdot)$ on $\R\times (0,1)$ from Proposition \ref{Wcont} that 
$$
\lim_{\mathcal{C}' \ni (z,\pi) \to (z_0,\pi_0) \in \partial\mathcal{C}'} \Big|\frac{\partial W}{\partial z}(z,\pi) \Big| < \infty 
\qquad \Leftrightarrow \qquad 
\frac{\partial W}{\partial z} \in L^\infty_{loc}(\R \times (0,1)).
$$
\vspace{0.15cm}
 
{\it Step 3: $\frac{\partial^2 w}{\partial \pi^2}(\cdot,\cdot)$ is bounded on the closure of $\mathcal{B} \cap \mathcal{C}'$, for all bounded sets $\mathcal{B}$.} 
Recall from Proposition \ref{prop3}.(iii), that $w$ solves the PDE \eqref{PDEw} on $\mathcal{C}'$, which implies in view of the definition \eqref{genZPi} that 
\begin{equation*}
\frac{1}{2}\Big(\frac{2\alpha}{\sigma}\Big)^2 \pi^2 (1-\pi)^2 \frac{\partial^2 w}{\partial \pi^2}(z,\pi)
=rw(z,\pi) - \frac{1}{2}\sigma^2\frac{\partial w}{\partial z}(z,\pi) - q(z,\pi) , \quad \text{for} \quad (z,\pi)\in\mathcal{C}'. 
\end{equation*}
Due to {\it step 2} and the smooth expression \eqref{W} of $F$, we observe that the right-hand side of the above expression is bounded on the closure of $\mathcal{B} \cap \mathcal{C}'$, for all bounded sets $\mathcal{B}$. 
Hence, $\frac{\partial^2 w}{\partial \pi^2}(\cdot,\cdot)$ is bounded on the closure of $\mathcal{B} \cap \mathcal{C}'$.
\vspace{0.15cm}

{\it Step 4: $\frac{\partial w}{\partial \pi}(z_0,\pi_0) =0$, for all $(z_0,\pi_0)\in\partial\mathcal{C}'$ such that $\pi_0 = c(z_0)$.}
Take now $(z_0,\pi_0)\in\partial\mathcal{C}'$, so that $\pi_0 = c(z_0)$. 
Due to {\it step 3} we know that the left-derivative $\frac{\partial w}{\partial \pi}(z,\pi-)$ exists. 
It thus follows from Proposition \ref{prop3}.(ii) that $\frac{\partial w}{\partial \pi}(z_0,\pi_0-) \leq 0$, since $(z_0,\pi_0-)\in\mathcal{C}'$ thanks to \eqref{C'S'c}, while we know from the definition \eqref{newval3} of $w$ that $\frac{\partial w}{\partial \pi}(z_0,\pi_0+) = 0$, since $(z_0,\pi_0+)\in\mathcal{S}'$.
Aiming for a contradiction we assume that $\frac{\partial w}{\partial \pi}(z_0,\pi_0-) < -\delta_0$, for some $\delta_0>0$. 

Then, we take a rectangular domain $\mathcal{B}$ around $(z_0,\pi_0)$ and define the stopping time
$$
\tau_\mathcal{B} := \inf\lbrace t\geq 0 \,|\, (Z^{z_0}_t,\Pi^{\pi_0}_t) \not\in \mathcal{B}\rbrace .
$$
In view of \eqref{Umart} together with \eqref{WU}, we know that 
\begin{equation*}
t\mapsto e^{-rt} W({Z}_{t},\Pi_{t}) \quad 
\text{is a supermartingale,}
\end{equation*}
which yields that 
\begin{align*} 
w(z_0,\pi_0) \geq \mathbb{E}_{z_0,\pi_0} \bigg[ e^{-r ({t\wedge \tau_\mathcal{B}})} w(Z_{t\wedge \tau_\mathcal{B}}, \Pi_{t\wedge \tau_\mathcal{B}}) + \int_0^{t\wedge \tau_\mathcal{B}} e^{-rs} q(Z_s, \Pi_s) ds \bigg]. 
\end{align*}
Given that $t \mapsto Z_{t\wedge \tau_\mathcal{B}}$ is increasing, $z \mapsto w(z,\pi)$ is non-decreasing on $\R$ due to Proposition \ref{prop3}.(i) and $q(\cdot,\cdot)$ is bounded on $\mathcal{B}$ by a constant 
$
c_\mathcal{B} := \sup_{(z,\pi)\in\mathcal{B}} \left| q(z,\pi) \right| \geq 0,
$
we have 
\begin{align*} 
w(z_0,\pi_0) \geq \mathbb{E}_{z_0,\pi_0} \bigg[ e^{-r ({t\wedge \tau_\mathcal{B}})} w(z_0, \Pi_{t\wedge \tau_\mathcal{B}}) - c_\mathcal{B} (t\wedge \tau_\mathcal{B}) \bigg]. 
\end{align*}
Then, we can use Tanaka's formula on $(e^{-r s} w(z_0, \Pi_s))_{s\in[0, t\wedge \tau_\mathcal{B}]}$ thanks to 
{\it step 3} to get 
\begin{align*} 
0 \geq \mathbb{E}_{z_0,\pi_0} \bigg[ \int_0^{t\wedge \tau_\mathcal{B}} \hspace{-3mm}e^{-rs}(\mathcal{G}-r)w(z_0, \Pi_s) {\bf 1}_{\{\Pi_s \not= \pi_0\}} ds 
+ \int_0^{t\wedge \tau_\mathcal{B}} \hspace{-3mm} e^{-rs}\big(w(z_0, \pi_0+) - w(z_0, \pi_0-)\big) dL_s^{\pi_0} 
- c_\mathcal{B} (t\wedge \tau_\mathcal{B}) \bigg],
\end{align*}
where $L^{\pi_0}$ is the local time of $\Pi$ at $\pi_0$. 
However, given that $\frac{\partial w}{\partial \pi}(z_0,\pi_0+) = 0$ and the assumption $\frac{\partial w}{\partial \pi}(z_0,\pi_0-) < -\delta_0$, as well as the boundedness of $(\mathcal{G}-r)w(\cdot, \cdot)$ on the closure of $\mathcal{B} \cap \mathcal{C}'$, we obtain for another constant $\bar c_\mathcal{B} \geq 0$ that
\begin{align*} 
0 > \mathbb{E}_{z_0,\pi_0} \bigg[ \delta_0 \int_0^{t\wedge \tau_\mathcal{B}} \hspace{-3mm} e^{-rs} dL_s^{\pi_0} 
- \bar c_\mathcal{B} (t\wedge \tau_\mathcal{B}) \bigg]
\geq \delta_0 e^{-rt} \,\mathbb{E}_{z_0,\pi_0} \left[ L_{t\wedge \tau_\mathcal{B}}^{\pi_0} \right] 
- \bar c_\mathcal{B} \,\mathbb{E}_{z_0,\pi_0} \left[t\wedge \tau_\mathcal{B} \right].
\end{align*}
This implies that 
\begin{align*} 
\delta_0 e^{-rt} \,\mathbb{E}_{z_0,\pi_0} \left[ L_{t\wedge \tau_\mathcal{B}}^{\pi_0} \right] 
< \bar c_\mathcal{B} \,\mathbb{E}_{z_0,\pi_0} \left[t\wedge \tau_\mathcal{B} \right],
\end{align*}
which leads to a contradiction for small enough $t$, since we can show that $\mathbb{E}_{z_0,\pi_0} [ L_{t\wedge \tau_\mathcal{B}}^{\pi_0} ] \sim \sqrt{t\wedge \tau_\mathcal{B}}$ by arguments similar to those in Lemma 13 of \cite{Peskir2019}.
\hfill \cvd

{\bf Proof of Theorem \ref{intw}.}
We will prove the two equalities sequentially. 

{\it Proof of $1^{st}$ equality.} 
Take $T>0$ and $(z,\pi)\in\R\times(0,1)$. 
Then, it follows from \cite[Theorem 3.1]{Peskir2005} that
\begin{align*}
e^{-r ({T \wedge \tau_n})} w(Z^z_{T\wedge \tau_n}, \Pi^\pi_{T\wedge \tau_n}) 
&= w(z,\pi) 
- \int_0^{T\wedge \tau_n} \hspace{-3mm}e^{-rs} q(Z^z_s, \Pi^\pi_s) {\bf 1}_{\{\Pi^\pi_s < c(Z^z_s)\}} ds \\&\quad 
+ \int_0^{T\wedge \tau_n}  \hspace{-3mm}e^{-rs} \frac{\partial w}{\partial \pi}(Z^z_s, \Pi^\pi_s) \Pi^\pi_s (1-\Pi^\pi_s) dW_s , 
\end{align*}
where we define for all $n\in\mathbb{N}$, the stopping times
$$
\tau_n := \inf \Big\{t\geq 0 \,\Big|\, \int_0^t e^{-rs} \Big(\frac{\partial w}{\partial \pi}(Z^z_s, \Pi^\pi_s) \Pi^\pi_s (1-\Pi^\pi_s)\Big)^2 ds \geq n \Big\}.
$$
Then, by taking expectations, we have 
\begin{align*} 
w(z,\pi) 
&= \mathbb{E}_{z,\pi} \bigg[e^{-r ({T \wedge \tau_n})} w(Z^z_{T\wedge \tau_n}, \Pi^\pi_{T\wedge \tau_n}) 
+ \int_0^{T\wedge \tau_n} \hspace{-3mm}e^{-rs} q(Z^z_s, \Pi^\pi_s) {\bf 1}_{\{\Pi^\pi_s < c(Z^z_s)\}} ds \bigg]. 
\end{align*}
Then, taking limits as $n \uparrow \infty$ and $T \uparrow \infty$, it follows from the dominated convergence theorem that
\begin{align} 
\label{eq:domconv}
w(z,\pi) 
&= \mathbb{E}_{z,\pi} \bigg[ \int_0^{\infty} e^{-rs} q(Z^z_s, \Pi^\pi_s) {\bf 1}_{\{\Pi^\pi_s \leq c(Z^z_s)\}} ds \bigg] , 
\end{align}
where the replacement of ${\bf 1}_{\{\Pi^\pi_s < c(Z^z_s)\}}$ with ${\bf 1}_{\{\Pi^\pi_s \leq c(Z^z_s)\}}$ is possible because $(\Pi_t^\pi)_{t\geq 0}$ admits an absolutely continuous transition density $(p_t(\pi,\pi'))_{t\geq 0, (\pi,\pi')\in(0,1)^2}$ due to \cite[Theorem 2.3.1]{Nualart} and $(Z_t^z)_{t\geq 0}$ is a deterministic process. 

{\it Proof of $2^{nd}$ equality.} 
This follows by expressing the expectation as an integral with respect to the transition density of $(\Pi_t^\pi)_{t\geq 0}$.
\hfill \cvd

\section*{Acknowledgments}

\noindent The second author acknowledges the funding by the Deutsche Forschungsgemeinschaft (DFG, German Research Foundation) – Project-ID 317210226 – SFB 1283.


\begin{thebibliography}{99}
\addcontentsline{toc}{section}{References}

\bibitem{Allenetal} Allen, M.R., Frame, D.J., Huntingford, C., Jones, C.D., Lowe, J.A., Meinshausen, M., Meinshausen, N.\ (2009). Warming caused by cumulative carbon emissions towards the trillionth tonne. \emph{Nature} 458, 1163-1166.

\bibitem{AthanassoglouXepapadeas} 
Athanassoglou, S., Xepapadeas, A.\ (2012). Pollution control with uncertain stock dynamics: When, and how, to be precautious. \emph{Journal of Environmental Economics and Management} 63(3), 304--320.

\bibitem{Bain} Bain, A., Crisan D.\ (2009). \textit{Fundamentals of stochastic filtering. Stochastic modelling and applied probability}, 60. Springer-Verlag, New-York.

\bibitem{Barnett} Barnett, M.\ (2022). Climate change and uncertainty: An asset pricing perspective. Forthcoming on \emph{Management Science}.

\bibitem{Borodin} Borodin, A.N., Salminen, P.\ (2015). \textit{Handbook of Brownian motion - Facts and Formulae} (2nd edition). Birkh\"auser.

\bibitem{Boucekkine-etal} Boucekkine, R., Fabbri, G., Federico, S., Gozzi, F.\ (2022). A dynamic theory of spatial externalities. \emph{Games and Economic Behavior} 132, 133--165.

\bibitem{CDP} Carbon Disclosure Project (CDP) (2023). Scoping out: Tracking nature across the supply chain. \emph{Global Supply Chain Report 2022}, March 2023.

\bibitem{DeAngelis1} De Angelis, T., Gensbittel, F., Villeneuve, S.\ (2021). A Dynkin game on assets with incomplete information on the return. \emph{Mathematics of Operations Research} 46(1), 28-60.

\bibitem{DeAngelis2} De Angelis, T.\ (2020). Optimal dividends with partial information and stopping of a degenerate reflecting diffusion. \emph{Finance and Stochastics} 24(1), 71-123.

\bibitem{Decamps1} D\'ecamps, J.P., Mariotti, T., Villeneuve, S.\ (2005). Investment timing under incomplete information. \emph{Mathematics of Operations Research} 30(2), 472-500.

\bibitem{Dalbyetal} Dalby, P.A.O., Gillerhaugen, G.R., Hagspiel, V., Leth-Olsen, T., Thijssen, J.J.J.\ (2018). Green investment under policy uncertainty and Bayesian learning. \emph{Energy} 161, 1262-1281.

\bibitem{DiPi} Dixit, A.K., Pindyck, R.S.\ (1994). \textit{Investment under uncertainty}. Princeton University Press (Princeton).

\bibitem{Embrechts} Embrechts, P., Hofert, M.\ (2013). A note on generalized inverses. \emph{Mathematical Methods of Operations Research} 77, 423-432.

\bibitem{Falboetal} Falbo, P., Ferrari, G., Rizzini, G., Schmeck, M.D.\ (2021). Optimal switch from a fossil-fueled to an electric vehicle. \emph{Decisions in Economics and Finance} 44, 1147-1178.

\bibitem{FFR} Federico, S.,  Ferrari, G., Rodosthenous, N.\ (2021). Two-sided singular control of an inventory with unknown demand trend.  \emph{SIAM Journal on Control and Optimization} 61(5), 3076-3101.

\bibitem{Ferrarietal} Ferrari, G., Li, H., Riedel, F.\ (2022). A Knightian Irreversible Investment Problem. \emph{Journal of Mathematical Analysis and Applications} 507(1), 125744.

\bibitem{FloraTankov} Flora, M., Tankov, P.\ (2023). Green investment and asset stranding under transition scenario uncertainty. \emph{Energy Economics}, 106773. 

\bibitem{HellmannThijssen} Hellmann, T., Thijssen, J.J.J.\ (2018). Fear of the market or fear of the competitor? Ambiguity in a real options game. \emph{Operations Research} 66(6), 1744-1759.

\bibitem{HLG21} Huang, W., Liang, J., Guo, H.\ (2021). Optimal investment timing for carbon emission reduction technology with a jump-diffusion process. \emph{SIAM Journal on Control and Optimization} 59(5), 4024-4050.

\bibitem{Jeanblanc} Jeanblanc, M., Yor, M., Chesney, M.\ (2009). \textit{Mathematical methods for financial markets}. Springer.

\bibitem{JP17} {Johnson, P., Peskir, G.}\ (2017). Quickest detection problems for Bessel processes. \emph{Annals of Applied Probability} {27}, 1003-1056.

\bibitem{Karatzas1} Karatzas, I., Shreve, S.E.\ (1998). \textit{Methods of mathematical finance}. Applications of Mathematics (New York), 39. Springer-Verlag, New York.

\bibitem{Karatzas2} Karatzas, I., Shreve, S.E.\ (1991). \textit{Brownian motion and stochastic calculus} (Second Edition). Graduate Texts in Mathematics 113, Springer-Verlag, New York.

\bibitem{Kim} Kim, S.H.\ (2015). Time to come clean? Disclosure and inspection policies for green production. \emph{Operations Research} 63(1), 1-20. 

\bibitem{Kwon1} Kim, Y., Kwon, H.D.\ (2022). Investment in the common good: Free rider effect and the stability of mixed strategy equilibria. Article in Advance on \emph{Operations Research}, \url{https://doi.org/10.1287/opre.2022.2371}

\bibitem{Kwon2} Kwon, H.D.\ (2022). Game of variable contributions to the common good under uncertainty. \emph{Operations Research} 70(3), 1359–1370.

\bibitem{Krylov} Krylov, N.V.\ (2008). \textit{Lectures on elliptic and parabolic equations in Sobolev spaces}. Graduate studies in Mathematics, AMS Providence.

\bibitem{Lappi} Lappi, P.\ (2018). Optimal clean-up of polluted sites. \emph{Resource and Energy Economics} 54, 53-68.

\bibitem{LaRiviere-etal} La Riviere, J., Kling, D., Sanchirico, J.N., Sims, C.,  Springborn, M.\ (2018). The treatment of uncertainty and learning in the economics of natural resource and environmental management. \emph{Review of Environmental Economics and Policy} 12(1).

\bibitem{Liptser} Liptser, R., Shiryaev, A.N.\ (2001). \textit{Statistics of random processes I. General Theory}. Springer-Verlag, Berlin Heidelberg.

\bibitem{McD} McDonald, R.\ and Siegel, D.\ (1986). The value of waiting to invest. \emph{The Quarterly Journal of Economics} 101(4), 707-728.

\bibitem{Murto} Murto, P.\ (2007). Timing of investment under technological and revenue-related uncertainties. \emph{Journal of Economic Dynamics and Control} 31(5), 1473--1497. 668-

\bibitem{NishimuraOzaki} Nishimura, K.G, Ozaki, H.\ (2007). Irreversible investment and Knightian uncertainty. \emph{Journal of Economic Theory} 136(1), 668-694.

\bibitem{Nordhaus2} Nordhaus, W.D.\ (1991). To slow or not to slow: The economics of the greenhouse effect. \textit{The Economic Journal}, 101, 920-937.

\bibitem{Nordhaus1} Nordhaus, W.D.\ (2007). A review of the Stern review on the economics of climate Change. \textit{Journal of Economic Literature}. \textbf{XLV} 686-702.

\bibitem{Nualart} Nualart, D.\ (2006). \textit{The Malliavin calculus and related topics}, Springer, $2^{nd}$ Edition.

\bibitem{Peskir2005} Peskir, G.\ (2005). A change-of-variable formula with local time on curves. \emph{Journal of Theoretical Probability} 18, 499-535.

\bibitem{PeskirAm} Peskir, G.\ (2005). On the American option problem. \emph{Mathematical Finance} 15, 169-181.

\bibitem{Peskir2019} Peskir, G.\ (2019). Continuity of the optimal stopping boundary for two-dimensional diffusions. \emph{Annals of Applied Probability} 29, 505-530.

\bibitem{Peskir1} Peskir, G., Shiryaev, A.N.\ (2006). \textit{Optimal stopping and free-boundary problems}. Lectures in Mathematics ETH, Birkhauser.

\bibitem{Pindyck1} Pindyck, R.S.\ (2000). Irreversibilities and the timing of environmental policy. \emph{Resource and Energy Economics} 22, 233-259.

\bibitem{Pindyck0} Pindyck, R.S.\ (2002). Optimal timing problems in environmental economics. \emph{Journal of Economic Dynamics and Control} 26(9--10), 1677-1697.

\bibitem{Pindyck2} Pindyck, R.S.\ (2007). Uncertainty in environmental economics. \emph{Review of Environmental Economics and Policy} 1, 45-65.

\bibitem{SBTi} Science Based Targets initiative (SBTi) (2023). Catalyzing value chain decarbonization. \emph{Corporate Survey Results}, February 2023.

\bibitem{Shi} Shiryaev, A.N.\ (1978), \textit{Optimal stopping rules}. Springer (New York–Heidelberg).

\bibitem{Shiryaev1} Shiryaev, A.N.(2010). Quickest detection problems: Fifty years later. \emph{Sequential Analysis} 29(4), 345-385.

\bibitem{Stern} Stern, N.\ (2006). \textit{The Economics of climate change: The Stern review}. Cambridge: Cambridge University Press.

\bibitem{BT} Thijssen, J.J.J., Bregantini,D.\ (2017). Costly sequential experimentation and project valuation with an application to health technology assessment. \emph{Journal of Economic Dynamics and Control} 77, 202--229.

\bibitem{Tol1} Tol, R.S.J.\ (2018). Economic impacts of climate change. \emph{Review of Environmental Economics and Policy} 12, 4-25.

\bibitem{Tol2} Tol, R.S.J.\ (2006). The Stern review of the economics of climate change: A comment. \emph{Energy \& Environment} 17, 977-981.

\bibitem{Weitzman} Weitzman, M.L.\ (2007). A review of the Stern review on the economics of climate change. \emph{Journal of Economic Literature} XLV, 703-724. 

\end{thebibliography}
\end{document}